%Authors: N.J. Kalton and M.I. Ostrovskii

%Title: Distances between Banach spaces

%Filename: kaltonostrovdist.tex
%TeX: AMSTeX
%Length: 79013 bytes
%Received Date: 9/8/97
%SubjectClass: 46B20
%Abstract: The  main  object  of  the  paper  is  to study the distance
%between  Banach  spaces  introduced  by  Kadets.  For Banach
%spaces $X$ and $Y$, the Kadets distance is defined to be the
%infimum of the  Hausdorff distance $d(B_X,B_Y)$  between the
%respective  closed  unit  balls  over  all  isometric linear
%embeddings of $X$  and $Y$ into  a common Banach  space $Z.$
%This is compared with the Gromov-Hausdorff distance which is
%defined to be the infimum of $d(B_X,B_Y)$ over all isometric
%embeddings  into  a  common  metric  space  $Z$.    We prove
%continuity type results for the Kadets distance including  a
%result  that  shows  that   this  notion  of  distance   has
%applications to the theory of complex interpolation.

%Citation: preprint

%Special character check block
%32   space        33 ! exclam. pt.   34 " double quote  35 # sharp
%36 $ dollar       37 % percent       38 & ampersand     39 ' prime
%40 ( left paren.  41 ) rt. paren.    42 * asterisk      43 + plus
%44 , comma        45 - minus         46 . period        47 / divide
%58 : colon        59 ; semi-colon    60 < less than     61 = equal
%62 > greater than 63 ? question mark 64 @ at
%91 [ left bracket 92 \ backslash     93 ] right bracket 94 ^ caret
%95 _ underline    96 ` left single quote
%123 { left brace  124 | vertical bar 125 } right brace  126 ~ tilda

%Insert your TeX file starting here.

\input amstex
\documentstyle{amsppt}
\NoBlackBoxes
\def\sgn{\text{sgn }}
\magnification=1200
\topmatter
\title Distances between Banach spaces
\endtitle
\rightheadtext{distances between Banach spaces}
\author N.J.Kalton and M.I.Ostrovskii
\endauthor
\address Department of Mathematics, University of Missouri-Columbia,
Columbia, Missouri 65211, USA\endaddress
\email nigel\@math.missouri.edu\endemail
\address Mathematical Division, Institute for Low Temperature
Physics, 47 Lenin avenue, 310164 Kharkov, UKRAINE
\endaddress
\email mostrovskii\@ilt.kharkov.ua
\endemail
\keywords Banach Space, Opening between Subspaces, Gap between
Subspaces, Complex Interpolation Family, Gromov-Hausdorff Distance,
Quasi-Banach Space \endkeywords
\subjclass Primary  46B20, 46M35, Secondary 46B03, 54E35
\endsubjclass

\abstract  The main object of the paper is to study the distance between
Banach spaces introduced by Kadets.
 For Banach spaces $X$ and $Y$, the Kadets distance is defined to be
the infimum of the Hausdorff distance $d(B_X,B_Y)$ between the
respective closed unit balls over all isometric
linear embeddings of $X$ and $Y$ into a common Banach space $Z.$  This is
compared with the Gromov-Hausdorff distance which is defined to be the
infimum of $d(B_X,B_Y)$ over all isometric embeddings into a common
metric space $Z$.
We prove continuity type results for the Kadets distance including a
result that shows that this notion of distance has applications to the
theory of complex interpolation.
\endabstract
\endtopmatter

\heading 1. Introduction   \endheading

The standard notion of distance between two Banach spaces is the {\it
Banach-Mazur distance} which is defined by
$$ d_{BM}(X,Y)=\log \inf\{ \|T\|\|T^{-1}\|:\  T:X\to Y \text{ is an
isomorphism}\}.$$
(It is usual to omit the logarithm, but for consistency we will include
it).  The Banach-Mazur distance is only finite when $X$ and $Y$ are
isomorphic.  The main object of this paper is to study a measure of
distance we call the  Kadets distance and certain related
notions of distance.  The Kadets distance has natural applications
in interpolation theory  which we explain.

We recall that if $Z$ is a Banach space and $X$ and $Y$ are
closed subspaces of
$Z$ the {\it gap} or {\it opening} $\Lambda(X,Y)$ is defined as the
Hausdorff distance between the closed unit balls $B_X$ and $B_Y$
of $X$ and $Y$ i.e.
$$ \Lambda(X,Y) =\max\{\sup_{y\in B_Y}d(y,B_X),\sup_{x\in
B_X}d(x,B_Y)\}.$$
If $X$ and $Y$ are arbitrary Banach spaces we define the {\it Kadets
distance}
$$ d_K(X,Y)=\inf_{Z,U,V}\Lambda(UX,VY)$$ where the infimum is taken over
all Banach spaces  $Z$ and all linear isometric embeddings $U:X\to Z$ and
$V:Y\to Z.$

This distance was apparently introduced by Kadets \cite{15} who proved
for
example that $\lim_{p\to 2}d_K(\ell_p,\ell_2)=0.$  However the basic idea
seems to be implicit in some earlier work  of Krein, Krasnoselskii
and Milman \cite{20}, Brown \cite{5} and Douady \cite{10}.  The second
author studied the notion in \cite{22} and proved that $d_K$ satisfies
the triangle law but that there are non-isomorphic Banach spaces $X$
and $Y$ for
which $d_K(X,Y)=0$ (thus $d_K$ is a ``pseudo-metric'').  In the same
paper there is a completeness result: if $(X_n)$ is a sequence of Banach
spaces Cauchy with respect to $d_K$ then there is a Banach space $X$ so
that $\lim_{n\to\infty}d_K(X_n,X)=0.$

There is a series of papers studying the general problem of the
identifying properties which are stable under small perturbations in the
Kadets distance (see \cite{2}, \cite{5} \cite{10}, \cite{20},
\cite{22}, \cite{23}, \cite{24} and \cite{25}).
Precisely a property
$\Cal P$ is called
{\it stable} if there exists $\epsilon>0$ so that if $X$ has $\Cal P$ and
$d_K(X,Y)<\epsilon$ then $Y$ has $\Cal P.$  We refer to the survey
article \cite{24} Chapter 6: a partial list of stable properties
includes reflexivity, super-reflexivity, B-convexity (nontrivial
Rademacher type), the Banach-Saks property, the alternate-signs Banach
Saks property and the property of not containing $\ell_1.$

The Kadets distance is clearly related to the notion of Gromov-Hausdorff
distance between metric spaces (see \cite{11},\cite{27}; the precise
definition is given in Section 2). It is natural to introduce the
Gromov-Hausdorff distance  between two Banach spaces $X$ and $Y$ as
$$ d_{GH}(X,Y) =\inf_{Z,U,V}d(UB_X,VB_Y)$$
where the infimum is taken over all isometries of $U:X\to Z$ and $V:Y\to
Z
$ into a common metric space $Z$ (here $d(UB_X,VB_Y)$ is the Hausdorff
distance between $UB_X$ and $VB_Y.$)
Thus the Gromov-Hausdorff distance is simply the nonlinear analogue of
the Kadets distance.  It is not difficult to see that this definition
coincides with computing the standard Gromov-Hausdorff distance between
the unit balls $B_X,B_Y$ as metric spaces. Let us remark at this point
that a related global notion was considered in \cite{4} and \cite{12}.

In this paper we first compare these two notions of distance. It is worth
pointing out that one must distinguish between the case of complex
scalars and real scalars, because there are examples
(\cite{3},\cite{17},\cite{31}) of complex Banach spaces which are
real-isometric but not even complex-isomorphic.   Gromov-Hausdorff
distance cannot distinguish complex structures.

We show that (for real scalars) while Gromov-Hausdorff distance is not
equivalent to
the Kadets distance, the two notions are equivalent if one restricts to
Banach spaces which are nice enough.
For example if $X$ is B-convex (i.e. has non-trivial Rademacher type) or
if $X^*$ embeds into an $\Cal L_1-$space then $d_{GH}(X_n,X)\to 0$
implies
$d_K(X_n,X)\to 0.$  If $X$ is isomorphic to either $c_0$ or
$\ell_{\infty}$ then $d_{GH}(X_n,X)\to 0$ implies $d_{BM}(X_n,X)\to 0.$
On the other hand $d_{GH}(\ell_p,\ell_1)\to 0$ as $p\to 1$ while
$d_{K}(\ell_p,\ell_1)=1$ if $p>1.$  The precise identification of the
class on which the two distances are equivalent is related to the notion
of a $\Cal K-$space introduced in \cite{16} (see \cite{18}) and to
the theory of twisted sums.

In fact for real scalars, Gromov-Hausdorff distance is equivalent to a
notion of distance analogous to the Kadets distance but allowing the
superspace $Z$ to be a quasi-Banach space.

These results are developed in Section 3, after some preliminary results
in Section 2.  In Section 4, we then apply our techniques to prove a
number of continuity-type results for the Kadets metric.  For example we
show that in an obvious sense the map $X\to X^*$ is continuous for the
Kadets metric, and even more one has $d_K(X^*,Y^*)\le 2d_K(X,Y)$ for any
pair of Banach spaces.  We also show that if $(X_0,X_1)$ is a
complex Banach couple and $X_{\theta}=[X_0,X_1]_{\theta}$ are the
spaces obtained by the (Caldero\'n) method of complex interpolation (cf.
\cite{6}) then the map $\theta\to X_{\theta}$ is continuous for the
Kadets distance for $0<\theta<1$. This result is closely related to
recent work on uniform homeomorphisms between the unit balls of two
Banach spaces using complex interpolation methods (cf. \cite{8}). We
give precise estimates here and obtain the estimate for $1<p,q<\infty$,
$$ d_K(\ell_p,\ell_q)\le
2\frac{\sin(\pi |1/p-1/q|/2)}{\sin(\pi (1/p+1/q)/2)} $$
which improves earlier estimates (see \cite {15}, \cite{22}).  We
remark
that in \cite{23} or \cite{24} (pp. 292, 303) there is a lower estimate
$d_K(\ell_p,\ell_q)\ge \frac12(2^{1/p}-2^{1/q}).$

Finally in Section 5, we make some remarks on the topology of the
pseudo-metric space of all Banach spaces with a given density character
with either notion of distance.  We point out that results on stability
or openness of some property lead automatically to results on complex
interpolation spaces, and also show that the continuity results of the
previous section lead to new stable or open properties.  We identify the
component of $\ell_1$ for the Kadets distance and raise the question of
identifying the components of $\ell_2$ and $c_0.$  We do not know if the
set of separable Banach spaces is connected for the Gromov-Hausdorff
distance.  We also show that if $1<p\neq 2<\infty$ the set of
spaces
isomorphic to $\ell_p$ is non-separable for both notions of distance.

\heading
2. Gromov-Hausdorff distance and Kadets distance.
\endheading

We first recall the notion of Gromov-Hausdorff distance between metric
spaces. It will be convenient to expand the definition to include
pseudo-metric spaces.  We recall that if $M$ is a set, a pseudo-metric on
$M$ is a map $d:M\times M\to [0,\infty)$ which is symmetric and satisfies
the triangle law, the condition $d(x,x)=0$,
but not necessarily the condition $d(x,y)=0$ implies
$x=y.$  Suppose
$A$ and
$B$ are metric spaces
(or pseudo-metric
spaces) with bounded metrics.  We define the Gromov-Hausdorff distance
between
$A$ and
$B$
denoted $d_{gh}(A,B)$ to be the infimum of all $\epsilon\ge 0$ so that
there is a pseudo-metric space $M$ and isometric embeddings $i_A:A\to
M$ and $i_B:B\to M$ such that the Hausdorff distance $d_M(i_AA,i_BB)\le
\epsilon.$  If $\Cal M$ is any set of metric spaces then $d_{gh}$
defines a pseudo-metric on $\Cal M$; note that $d_{gh}(A,B)=0$ does not
necessarily imply that $A$ and $B$ are isometric (unless they
are compact).

We will be interested in an alternative formulation.  For convenience
we
denote by $\Cal F(A,B)$ the collection of all pairs $(\phi,\psi)$ of
maps $\phi:A\to B$ and $\psi:B\to A.$
If $(\psi,\phi)$ is such pair let $G=G(\phi,\psi)$ be the union of the
graphs of $\psi,\phi.$  We define
$D(\phi,\psi)$ to be the supremum of all quantities
$\frac12|d_A(a_1,a_2)-d_B(b_1,b_2)|$ over all $a_1,a_2\in A,b_1,b_2\in
B$ with
$(a_i,b_i)\in G$ for $i=1,2.$

In the special case when $\phi$ is invertible and $\psi=\phi^{-1}$,
$D(\phi,\psi)$ is the supremum of all quantities $\frac12|d_B(\phi
a_1,\phi a_2)-d_A(a_1,a_2)|$ where $a_1,a_2\in A.$

\proclaim{Theorem 2.1}Let $A$ and $B$ be bounded metric spaces.  Then
$$d_{gh}(A,B)=\inf_{(\phi,\psi)\in\Cal
F}D(\phi,\psi).$$
\endproclaim

\demo{Proof}First suppose $M$ is a pseudo-metric space and that
$A$ and $B$ are isometrically embedded in $M.$
 Let
$\epsilon$ be the
Hausdorff distance between $A$ and $B.$  Then if $\sigma>\epsilon$
we can define $\phi:A\to B$ and $\psi:B\to A$ so that
$d_M(a,\phi a),d_M(\psi b,b)<\sigma$ for all $a\in A$ and
$b\in B.$  Now for $(a_i,b_i)\in G=G(\phi,\psi)$ we clearly have
$|d_A(a_1,a_2)-d_B(b_1,b_2)|<2\sigma.$

To obtain the converse direction suppose $(\phi,\psi)$ given and that
$D(\psi,\phi)=\sigma.$  We let $M=A\cup B$ (disjoint union) and
define
a pseudo-metric $d_M$ as follows.  We let $d_M$ coincide with $d_A$ on
$A$ and with $d_B$ on $B$.  If $a\in A$ and $b\in B$ then
$$ d_M(a,b)=\inf_{(a',b')\in G} \left( d_A(a,a')+d_B(b',b) \right)
+\sigma.$$

 We must check this is a pseudo-metric on $A\cup B$
$a_1,a_2\in A$ and $b\in B.$  We check that $d_A(a_1,a_2)\le
d_M(a_1,b)+d_M(a_2,b).$  To do this suppose that $(\alpha_1,\beta_1)$
and $(\alpha_2,\beta_2)\in G.$ Then
$$
\align
d_B(\beta_1,b)+d_B(\beta_2,b) &\ge d_B(\beta_1,\beta_2) \\
&\ge d_A(\alpha_1,\alpha_2)-2\sigma.\endalign $$
Hence
$$ d_A(a_1,\alpha_1)+d_B(\beta_1,b)+d_A(\alpha_2,a_2)+d_B(\beta_2,b)\ge
d_A(a_1,a_2)-2\sigma$$
which establishes our claim.

We can also show that $d_M(a_1,b)\le d_A(a_1,a_2)+d_M(a_2,b).$   We omit
the
details which are easy.  Arguing symmetrically with $A,B$ interchanged
gives the conclusion that $d_M$ is a pseudo-metric.  Clearly
$d_M(a,\phi a) \le \sigma$ and $d_M(\psi b,b)\le \sigma.$
This shows that $d_{gh}(A,B)\le \sigma.$\enddemo

If $X$ and $Y$ are Banach spaces we define the Gromov-Hausdorff distance
$d_{GH}(X,Y)$ to be the Gromov-Hausdorff distance between their closed
unit balls $B_X,B_Y,$ i.e. $d_{GH}(X,Y)=d_{gh}(B_X,B_Y).$ Equivalently
$d_{GH}(X,Y)$ is the infimum of the
Hausdorff distance $d(B_X,B_Y)$ over all isometric embeddings of $X,Y$
into a common metric space $M.$
To establish this last comment, suppose $d$ is any metric on the formal
union $B_X\cup B_Y$, which coincides with the respective norm-distances
on $B_X$
and $B_Y$.  We can extend $d$ to $X\cup Y$ by defining $d$ again to
coincide with the norm-distance on each of $X$ and $Y$ and for $x\in
X,y\in Y$,
$$ d(x,y)=\inf_{u\in B_X,v\in B_Y}\{\|x-u\|_X + d(u,v)+\|y-v\|_Y\}.$$
We leave the details to the reader.

Let us note here that our definition applies to both real and complex
Banach spaces, but there are complex Banach spaces which are
real-isometric and not even complex-isomorphic (\cite {3},
\cite{17}, \cite{31}).
In
view of this, Gromov-Hausdorff distance is most natural for the category
of real Banach spaces.

\proclaim{Corollary 2.2}If $X$ and $Y$ are Banach spaces and
$d_{GH}(X,Y)<\sigma$ then there exist maps $\phi:B_X\to B_Y$ and
$\psi:B_Y\to B_X$ such that:
$$ |\|x-\psi(y)\|_X-\|y-\phi(x)\|_Y|< 2\sigma$$
whenever $x\in B_X,y\in B_Y.$\endproclaim

The Kadets distance $d_K(X,Y)$ is defined to be the infimum of
the gap  $\Lambda(X,Y) (=d(B_X,B_Y))$ over all linear isometric
embeddings of
$X$ and
$Y$ into a
common Banach space $Z.$ Here our definition can equally be applied to
the real or complex case. The Kadets distance is again a pseudo-metric
on any set of Banach spaces (see \cite{15} \cite{22}).  We clearly
have
the inequality $d_{GH}(X,Y)\le d_K(X,Y).$

We now give a similar formulation of the
Kadets distance.  Let $\Cal F_h(X,Y)$ be the set of all
pairs of homogenous maps $\Phi:X\to Y$ and $\Psi:Y\to X$ such that
$\|\Phi(x)\|_Y\le \|x\|_X$ and $\|\Psi(y)\|_X\le \|y\|_Y$ for all $x\in
X$ and $y\in Y.$  We define $\Delta=\Delta(\Phi,\Psi)$ to be the
least constant such that
$$\left
|\|\sum_{i=1}^mx_i-\sum_{j=1}^n\Psi(y_j)\|_X-\|\sum_{i=1}^m\Phi(x_i)
-\sum_{j=1}^ny_j\|_Y \right|\le
\Delta\left(\sum_{i=1}^m\|x_i\|_X+\sum_{j=1}^n\|y_j\|_Y\right)$$ for all
$x_1,\ldots,x_m\in X$ and $y_1,\ldots, y_n\in Y.$

In the special case when $\Psi=\Phi^{-1}$ (and hence $\Phi$ is
norm-preserving) notice $\Delta$ is the least constant such that
$$ \left|\|\sum_{i=1}^m\Phi(x_i)\|_Y-\|\sum_{i=1}^mx_i\|_X\right|\le
\Delta\sum_{i=1}^m\|x_i\|_X.$$

\proclaim{Theorem 2.3}If $X,Y$ are Banach spaces then
$$d_K(X,Y)=\inf_{
(\Phi,\Psi)\in\Cal F_h} \Delta(\Psi,\Phi).$$
\endproclaim

\demo{Proof}If $X$ and $Y$ are isometrically embedded in $Z$ and
$d(B_X,B_Y)<\sigma$ we can construct $(\Phi,\Psi)\in\Cal F_h$ with
$\|x-\Phi(x)\|_Z\le \sigma \|x\|_X$ and $\|y-\Psi(y)\|_Y\le
\sigma\|y\|_Y$ whenever $x\in X$ and $y\in Y.$  It is then trivial that
$\Delta(\Phi,\Psi)\le \sigma.$

For the converse direction, let us suppose that $(\Phi,\Psi)\in\Cal
F_h(X,Y)$ are given and that $\Delta(\Phi,\Psi)=\sigma.$  We will define
$Z$ to be the direct sum $X\oplus Y$ equipped with an equivalent norm.
Precisely we define
$$ \|(u,v)\|_Z=\inf\left\{\|x_0\|_X +\|y_0\|_Y
+\sigma(\sum_{i=1}^m\|x_i\|_X+\sum_{j=1}^n\|y_j\|_Y)\right\}$$
where the infimum is taken over all $\{x_i\}_{i=0}^m$ in $X$ and
$\{y_j\}_{j=0}^n$ in $Y$ such that
$$ u=x_0+\sum_{i=1}^mx_i +\sum_{j=1}^n\Psi(y_j)$$
and
$$ v=y_0+\sum_{i=1}^m\Phi(x_i)+\sum_{j=1}^ny_j.$$

We show that $\|(u,0)\|_Z=\|u\|_X$  (and then it follows similarly that
$\|(0,v)\|_Z=\|v\|_Y.$)  Indeed  suppose
$$ u=x_0+\sum_{i=1}^mx_i+\sum_{j=1}^n\Psi(y_j)$$ and
$$ 0=y_0 +\sum_{i=1}^m\Phi(x_i) +\sum_{j=1}^ny_j.$$

It follows from our condition that
$$ \|\sum_{i=1}^mx_i+\sum_{j=1}^n\Psi(y_j)\|_X \le \|y_0\|_Y
+\sigma(\sum_{i=1}^m\|x_i\|_X+\sum_{j=1}^n\|y_j\|_Y).$$
Hence
$$ \|u\|_X \le \|x_0\|_X +\|y_0\|_Y
+\sigma(\sum_{i=1}^m\|x_i\|_X+\sum_{j=1}^n\|y_j\|_Y).$$
This shows that $\|(u,0)\|_Z=\|u\|_X.$

Clearly the construction yields that $\Lambda(X,Y)\le \sigma$ since
$\|(x,\Phi(x))\|_Z \le \sigma \|x\|_X$ and $\|(\Psi(y),y)\|_Z\le
\sigma\|y\|_Y.$\qed\enddemo

There is an amusing way to interpret $\Delta(\Phi,\Psi).$  Form the
vector sequence spaces $\ell_1(X)$ and $\ell_1(Y)$ and define
$\tilde\Phi((x_n))= (\Phi(x_n))$ and $\tilde\Psi((y_n))=(\Psi(y_n)).$
 We can consider the  pair
 $(\tilde \Phi,\tilde\Psi)\in\Cal F_h(\ell_1(X),\ell_1(Y)).$  Then
consider the unit ball $B_{\ell_1(X)}$ with the pseudo-metric induced by
the seminorm $|(x_n)|_X=\|\sum_{n=1}^{\infty}x_n\|_X$; similarly
consider $B_{\ell_1(Y)}$ with the seminorm $|(y_n)|_Y$ defined in the
analogous way.

Then $\Delta(\Phi,\Psi) =D(\tilde\Phi,\tilde\Psi)$ for the unit balls of
$\ell_1(X)$ and $\ell_1(Y)$ equipped with the pseudo-metrics induced by
these seminorms.  We leave the details to the reader.

It is natural to ask for a theorem of the type of Theorem 2.3 but with
$\Phi$ a bijection and $\Psi=\Phi^{-1}.$  This can be done at the cost
of a constant in the calculations.

\proclaim{Theorem 2.4}Let $X$ and $Y$ be a Banach spaces with
$d_K(X,Y)<\sigma<\frac16$.  Then there is a norm-preserving bijection
$\Omega:X\to Y$ such that if $x_1,\ldots,x_n\in X$ then
$$ \left|\|\sum_{i=1}^n\Omega(x_i)\|_Y-\|\sum_{i=1}^nx_i\|_X\right|\le
14\sigma\sum_{i=1}^n\|x_i\|_X.$$ \endproclaim

\demo{Proof}First we observe that if $d_K(X,Y)<1/2$ then
$X$ and
$Y$
have the same density character (cf. \cite{20} or \cite{24} 6.23).  We
will suppose that
$X$ and $Y$ are simultaneously embedded into a common Banach space
$Z$ with $d(B_X,B_Y)=\delta <\sigma.$ Let us pick a maximal collection of
vectors
$(x_i:i\in I)$ in $S_X$ such that if $|a|=1$ and $i\neq j$ then
$\|x_i-ax_j\|_X> 4\sigma$.  Then there exist vectors
$(y_i:i\in I)$ in
$S_Y$ with $\|x_i-y_i\|_Z< 2\delta.$  If $y\in S_Y$ then there
exists
$x\in S_X$ with $\|y-x\|_Z<2\delta$ and so there exists $i\in I$ and
$|a|=1$ such that
$\|y-ay_i\|_Z< 8\sigma.$  On the other hand if $i\neq j$
and
$|a|=1$ then $\|y_i-ay_j\|_Z> 4(\sigma-\delta).$

We can now partition $S_X$ into sets $(A_i:i\in I)$ such that if $x\in
A_i$ and $|a|=1$ then $ax\in A_i$ and further that $\|x-x_i\|_X\le
2\sigma$ implies $x\in A_i$ while $x\in A_i$ implies
that there exists $a$ with $|a|=1$ and $\|x-ax_i\|_X \le
4\sigma.$
We then define $A'_i$ by taking one representative $u$ from each set
$\{ax:|a|=1\}$ contained in $A_i$ with the property $\|u-x_i\|_X \le
4\sigma.$

In the same way we
 can  partition $S_Y$ into sets $(B_i:i\in I)$ such that if $y\in
B_i$ and $|a|=1$ then $ay\in B_i$ and further that $\|y-y_i\|_X\le
2(\sigma-\delta)$ implies $y\in B_i$ while $y\in
A_i$ implies that there exists $a$ with $|a|=1$ and $\|y-ay_i\|_X \le
8\sigma.$
We then define $B'_i$ by taking one representative $v$ from each set
$\{ay:|a|=1\}$ contained in $B_i$ with the property $\|v-y_i\|_X \le
8\sigma.$

It is easy to see that the sets $A'_i,B'_i$ have the same cardinality for
each $i$ and so we can define a bijection $\Omega:S_X\to S_Y$ such that
$\Omega(A'_i)=B'_i$, $\Omega(ax)=a\Omega(x)$ when $|a|=1$.  If
$x\in S_X$ then there exists $i\in I$ and $|a|=1$ with
$ax\in A_i'$.  Thus $\Omega(ax)\in B_i'$ and so $\|x-\Omega(x)\|_Z \le
14\sigma.$

The result now follows immediately.\qed\enddemo

Let us give one immediate application of Theorem 2.3.

\proclaim{Theorem 2.5}Let $X$ and $Y$ be real Banach spaces and suppose
$Z$ is a metric linear space equipped with an invariant metric $d_Z$.
Suppose $X$ and $Y$ are linearly and isometrically embedded into $Z.$
Then
$d_K(X,Y)\le d_Z(B_X,B_Y).$\endproclaim

\demo{Remark}Thus it would make no difference in the definition of
$d_K$ to allow $Z$ to be a metric linear space instead of a Banach
space.\enddemo

\demo{Remark}The proof given below can be extended to complex spaces, if
we further assume that the metric $d_Z$ is invariant under
multiplication by $e^{i\theta}$ for $0<\theta<2\pi.$  For the
real case invariance under multiplication by $-1$ follows from
additive invariance since $d_Z(-x,-y)=d_Z(x-y,0)=d_Z(x,y).$\enddemo

\demo{Proof}Suppose $\sigma>d_Z(B_X,B_Y)$.  Then we can define
$\Phi:S_X\to B_Y$ and $\Psi:S_Y\to B_X$ such that $d_Z(x,\Phi(x))\le
\sigma$ and $d_Z(\Psi(y),y)\le \sigma$ for all $x,y.$  Since the metric
is translation-invariant it is clear that we can suppose
$\Phi(-x)=-\Phi(x)$ and $\Phi(-y)=\Phi(y).$  We then extend $\Phi$ and
$\Psi$ to be homogenous.

Now if $x_1,\ldots,x_m\in S_X$ and $y_1,\ldots, y_n\in S_Y$ then
$$
d_Z(\sum_{i=1}^mx_i+\sum_{j=1}^n\Psi(y_j),\sum_{i=1}^m\Phi(x_i)+
\sum_{j=1}^n y_j) \le (m+n)\sigma.$$

It follows that since $\|x\|_X=d_Z(x,0)$ for $x\in X$ and
$\|y\|_Y=d_Z(y,0)$ for $y\in Y$ we have:
$$
\left|\|\sum_{i=1}^mx_i+\sum_{j=1}^n\Psi(y_j)\|_X-\|\sum_{i=1}^n
\Phi(x_i)+\sum_{j=1}^ny_j\|_Y\right|\le (m+n)\sigma.$$
From this it follows easily that if $r_1,\ldots,r_m,s_1,\ldots,s_n$ are
integers that
$$
\left|\|\sum_{i=1}^mr_ix_i+\sum_{j=1}^n\Psi(s_jy_j)\|_X-\|\sum_{i=1}^n
\Phi(r_ix_i)+\sum_{j=1}^ns_jy_j\|_Y\right|\le
(\sum_{i=1}^mr_i+\sum_{j=1}^ns_j)\sigma.$$
Clearly the same inequality then holds for
$r_1,\ldots,r_m,s_1,\ldots,s_n$ rational and then even real by a density
argument.  This implies that $\Delta(\Phi,\Psi)\le \sigma.$\qed\enddemo

Now suppose $0<r<1.$
We recall that an $r$-norm on a real or complex vector space $X$ is a map
$x\to \|x\|_X$ such that:\newline
(1) $\|x\|_X> 0$ if $x\neq 0,$\newline
(2) $\|\alpha x\|_X=|\alpha|\|x\|_X$ for $\alpha\in \Bbb K,\ x\in
X$\newline
(3) $\|x_1+x_2\|_X^r\le \|x_1\|_X^r+\|x_2\|_X^r$ for $x_1,x_2\in X.$
\newline
Here $\Bbb K=\Bbb R$ or $\Bbb K=\Bbb C.$
Notice that in \cite{18} this is called an $r$-subadditive quasi-norm.
If the metric $d(x_1,x_2)=\|x_1-x_2\|_X^r$ makes $X$ complete we say that
$X$ is an $r$-normed quasi-Banach space.

Now suppose $X$ and $Y$ are Banach spaces (so that $X$ and $Y$ are
also $r$-normed quasi-Banach spaces for any $0<r<1.$)
We define
$d_{r}(X,Y)$ to be the infimum of
$d(B_X,B_Y)$ over all linear isometric embeddings of $X$ and $Y$ into an
$r$-normed quasi-Banach space $Z.$  Note here that the $r$-norm does not
induce a metric on $Z$ so this is not covered by the preceding theorem.
There is, however an analogue of Theorem 2.3 for this situation.  If
$(\Phi,\Psi)\in \Cal F_h$ then we define $\Delta_r(\Phi,\Psi)$ to be the
least constant $\Delta_r$ such that we have for any
$x_1,\ldots,x_m\in X$ and $y_1,\ldots,y_n\in Y$ that
$$ \left|\|\sum_{i=1}^mx_i-\sum_{j=1}^n\Psi(y_j)\|_X^r-\|\sum_{i=1}^m
\Phi(x_i)-\sum_{j=1}^ny_j\|_Y^r\right|\le
\Delta^r_r(\sum_{i=1}^m\|x_i\|_X^r
+\sum_{j=1}^n\|y_j\|_Y^r).$$

In an exactly analogous fashion we may prove:

\proclaim{Theorem 2.6}If $X$ and $Y$ are Banach spaces then
$$d_{r}(X,Y)=\inf_{
(\Phi,\Psi)\in\Cal F_h} \Delta_r(\Psi,\Phi).$$  \endproclaim

\heading{3. Comparison of metrics} \endheading

\proclaim{Proposition 3.1}Suppose $X$ and $Y$ are Banach spaces and
$(\Phi,\Psi)\in \Cal F_h(X,Y).$  Let $\sigma=\Delta(\Phi,\Psi).$
Then \newline
(1) Given $y\in Y$ there exists $x\in X$ with $\|x\|_X\le \|y\|_Y$
and $\|y-\Phi(x)\|_Y\le 2\sigma\|y\|_Y$ \newline
(2) If $x\in X$ then $\|\Phi(x)\|_Y \ge (1-\sigma)\|x\|_X.$
\newline
(3) If $x_1,\ldots,x_n\in X$ and $\sum_{k=1}^nx_k=0$ then
$$\|\sum_{k=1}^n\Phi(x_k)\|_Y \le \sigma\sum_{k=1}^n\|x_k\|_X.$$
\endproclaim

\demo{Proof}(1) Just take $x=\Psi(y).$  Then $\|x-\Psi(y)\|_X=0$ so that
$\|y-\Phi(x)\|_Y \le \sigma(\|x\|_X+\|y\|_Y)\le 2\sigma\|y\|_Y.$
(2) and (3) are immediate from the definition of $\Delta(\Phi,\Psi).$
\qed\enddemo

In a very similar way we can establish:

\proclaim{Proposition 3.2}Suppose $0<r<1.$ Suppose $X$ and $Y$ are
Banach spaces and
$(\Phi,\Psi)\in \Cal F_h(X,Y).$  Let $\sigma=\Delta_r(\Phi,\Psi).$
Then \newline
(1.) Given $y\in Y$ there exists $x\in X$ with $\|x\|_X\le \|y\|_Y$
and $\|y-\Phi(x)\|_Y\le 2^{1/r}\sigma\|y\|_Y$ \newline
(2.) If $x\in X$ then $\|\Phi(x)\|_Y \ge (1-\sigma^r)^{1/r}\|x\|_X.$
\newline
(3.) If $x_1,\ldots,x_n\in X$ and $\sum_{k=1}^nx_k=0$ then
$$\|\sum_{k=1}^n\Phi(x_k)\|_Y \le
\sigma(\sum_{k=1}^n\|x_k\|^r_X)^{1/r}.$$
\endproclaim

Our interest in Propositions 3.1 and 3.2 is to establish a converse.

\proclaim{Proposition 3.3}Let $X$ and $Y$ be Banach spaces and suppose
$\Phi:X\to Y$ is a homogeneous map satisfying
$\frac12\|x\|_X \le \|\Phi(x)\|_Y\le
\|x\|_X$ such that for a constant $0\le \sigma< 1$ we have:\newline
(1) Given $y\in Y$ there exists $x\in X$ with $\|x\|_X\le \|y\|_Y$
and $\|y-\Phi(x)\|_Y\le \sigma\|y\|_Y$ \newline
(2) If $x_1,x_2,x_3\in X$ and $\sum_{k=1}^3x_k=0$ then
$$\|\sum_{k=1}^3\Phi(x_k)\|_Y \le \sigma\sum_{k=1}^3\|x_k\|_X.$$
\newline
Then if $\Psi:Y\to X$ is a homogenous map satisfying
$\|\Psi(y)\|_X\le \|y\|_Y$ and $\|y-\Phi(\Psi(y))\|_Y\le \sigma \|y\|_Y$
(whose existence is guaranteed by (1)) we have
$$ \left|\|x-\Psi(y)\|_X-\|y-\Phi(x)\|_Y\right|\le
6\sigma(\|x\|_X+\|y\|_Y).$$ Furthermore for each $0<r<1$ there is a
universal constant $C=C(r)$ such that $\Delta_r(\Phi,\Psi)\le
 C\sigma.$

If further we have
:\newline
(3) If $x_1,\ldots,x_n\in X$ and $\sum_{k=1}^nx_k=0$ then
$$\|\sum_{k=1}^n\Phi(x_k)\|_Y \le
\sigma\sum_{k=1}^n\|x_k\|_X,$$ then $\Delta(\Phi,\Psi)\le
20\sigma.$
\endproclaim

\demo{Proof}Suppose $x\in X$ and $y\in Y.$  Then
$$ \|\Phi(x)-\Phi(\Psi(y))-\Phi(x-\Psi(y))\|_Y \le
2\sigma(\|x\|_X+\|y\|_Y).$$
Hence
$$ \|\Phi(x)-y-\Phi(x-\Psi(y))\|_Y \le 3\sigma
(\|x\|_X+\|y\|_Y).$$
It follows that
$$ \|\Phi(x)-y\|_Y \le \|x-\Psi(y)\|_X + 3\sigma(\|x\|_X+\|y\|_Y)$$
and
$$ \|x-\Psi(y)\|_X \le \|\Phi(x)-y\|_Y +
6\sigma(\|x\|_X+\|y\|_Y).$$
This proves the first part of the Proposition.

Notice that the same proof yields that for $0<r<1,$
$$ \left|\|x-\Psi(y)\|^r_X-\|y-\Phi(x)\|^r_Y\right|\le
6^r\sigma^r(\|x\|_X+\|y\|_Y)^r.$$

Now suppose $0<r\le 1$ and that there is a constant $\tau$ so that for
any $x_1,\ldots, x_n\in X$ with $\sum_{i=1}^nx_i=0$ we have
$$ \|\sum_{i=1}^n\Phi(x_i)\|_Y \le
\tau(\sum_{i=1}^n\|x_i\|^r_X)^{1/r}.\tag 3.1$$

First suppose $y_1,\ldots,y_n\in Y$ and $\sum_{i=1}^ny_i=0.$
Let $v=\sum_{i=1}^n\Psi(y_i).$  Then
$$ \|\sum_{i=1}^n\Phi(\Psi(y_i))-\Phi(v)\|_Y \le 2^{1/r}\tau
(\sum_{i=1}^n\|y_i\|_Y^r)^{1/r}.$$

Thus
$$ \|\sum_{i=1}^ny_i-\Phi(v)\|_Y \le (2^{1/r}\tau+\sigma)
(\sum_{i=1}^n\|y_i\|_Y^r)^{1/r}.$$
It follows that $$\|v\|_X\le 2\|\Phi(v)\|_Y \le
(2^{1+1/r}\tau+2\sigma)(\sum_{i=1}^n\|y_i\|_Y^r)^{1/r}.$$

Now suppose that $x_1,\ldots,x_m\in X$ and $y_1,\ldots,y_n\in Y.$  Let
$u=\sum_{i=1}^mx_i$ and $v=\sum_{j=1}^ny_j.$  Then
$$ \left|\|u-\Psi(v)\|^r_X-\|v-\Phi(u)\|^r_Y\right|\le
6^r\sigma^r(\sum_{i=1}^m\|x_i\|_X^r
+\sum_{j=1}^n\|y_j\|_Y^r).$$
However we also have:
$$ \|\Psi(v)-\sum_{j=1}^n\Psi(y_j)\|_Y^r \le 2(2^{1+1/r}\tau+2\sigma)^r
\sum_{j=1}^n\|y_j\|_Y^r$$ and
$$ \|\Phi(u)- \sum_{i=1}^m\Phi(x_i)\|_X^r \le 2\tau^r
\sum_{i=1}^m\|x_i\|_X^r.$$

Combining these we have
$$ \Delta_r(\Phi,\Psi) \le (6^r\sigma^r
+2\tau^r+2(2^{1+1/r}\tau+2\sigma)^r)^{1/r}.$$

The final part of the Proposition is then immediate taking $r=1$ and
$\tau=\sigma.$

For the remaining case we observe that for
$0<r<1$ (2) implies that (3.1) holds with
a constant $\tau=C\sigma$ where $C$ depends on $r.$  This is a
well-known calculation first observed in \cite{16} (see also
\cite{18} p.91).  Suppose $x_1,\ldots,x_n$ are nonzero in $X$ and
$x^*\in X^*$ with $\|x^*\|_{X^*}=1$.  Define a linear map $T:\ell_r^n\to
X$ by $Te_i=x_i/\|x_i\|_X$;  then $\|T\|\le 1$ and
(2) implies that the functional
$F(\xi)=x^*(\Phi T(\xi))$ satisfies
$$ |F(\xi_1+\xi_2)-F(\xi_1)-F(\xi_2)|\le
2\sigma(\|\xi_1\|_r+\|\xi_2\|_r).$$  Appealing to Lemma 5.8 of
\cite{18}
we have
$$ |x^*(\Phi(\sum_{i=1}^nx_i))-\sum_{i=1}^n x^*(\Phi(x_i))|\le
2(\sum_{i=1}^n(2/i)^{1/r})\sigma(\sum_{i=1}^n\|x_i\|_X^r)^{1/r}.$$
This
establishes the Proposition, since
$\sum_{i=1}^{\infty}(2/i)^{1/r}< \infty.$\qed\enddemo

\proclaim{Lemma 3.4}Suppose $\epsilon>0.$  Then there exists
$\sigma=\sigma(\epsilon)>0$ so that if $X$ and $Y$ are real Banach spaces
and
$(\phi,\psi)\in\Cal F(B_X,B_Y)$ with $D(\phi,\psi)<\sigma$ then there exists
a homogeneous map $\Phi:X\to Y$ such that
$\frac12\|x\|_X\le \|\Phi(x)\|_Y\le
\|x\|_X$ for
$x\in X$ and:\newline
\item{1.} Given $y\in Y$ there exists $x\in X$ with $\|x\|_X\le \|y\|_Y$
and $\|y-\Phi(x)\|_Y\le \epsilon\|y\|_Y$ \newline
\item{2.} If $x_1,x_2,x_3\in X$ and $\sum_{k=1}^3x_k=0$ then
$$\|\sum_{k=1}^3\Phi(x_k)\|_Y \le \epsilon(\sum_{k=1}^3\|x_k\|_X).$$
\endproclaim

\demo{Proof}Suppose first $D(\phi,\psi)=\sigma.$  Let $\phi(0)=y\in B_Y.$
Then there exists $v\in B_Y$ with $\|y-v\|_Y=1+\|y\|_Y.$  Thus
$\|\psi(v)\|_X \ge 1+\|y\|_Y- 2\sigma.$  It follows that $\|\phi(0)\|_Y
\le 2\sigma.$  From this it follows similarly that for any $x\in B_X$ we
have $\|x\|_X-4\sigma\le\|\phi(x)\|_Y\le \|x\|_X+4\sigma.$

Now let us assume the conclusion of the Lemma is false.

Indeed if this is so we can find a sequence of pairs of Banach
spaces
$(X_n,Y_n)$ and pairs of functions $(\phi_n,\psi_n)$ with
$D(\phi_n,\psi_n)<\frac1{n^2}$ but such that for every homogeneous map
$\Phi:X_n\to Y_n$ with $\frac12\|x\|_{X_n}\le \|\Phi(x)\|_{X_n}\le
\|x\|_{X_n}$ either
(1) or
(2) fails.

We will in particular define $\Phi_n:X\to Y$ to be a homogenous map such
that either $\Phi_n(x)=(1+\frac4n)^{-1}\phi_n(x)$ or
$\Phi_n(x)=-(1+\frac4n)^{-1}\phi_n(-x)$ when
$\|x\|_X=\frac1n.$
As long as $n\ge 20$ we then have $\frac12\|x\|_{X_n}\le
\|\Phi_n(x)\|_{Y_n}\le \|x\|_{X_n}.$

Now let $\Cal U$ be a nonprincipal ultrafilter on $\Bbb N.$
  We form the
ultraproducts $X_{\Cal U}=\prod_{\Cal U}(X_n)$ and $Y_{\Cal
U}=\prod_{\Cal U}(Y_n).$  thus $X_{\Cal U}$ can be realized as the
Hausdorff quotient of the seminormed space $\ell_{\infty}(X_n)$ with the
seminorm $\|\bold x\|_{X_{\Cal U}}=\lim_{n\in\Cal U}\|x_n\|_{X_n}$ where
$\bold x=(x_n).$  We refer to \cite{13} or \cite{9} for details about
ultraproducts.

Define a map $\Omega:X_{\Cal U}\to Y_{\Cal U}$ by $$\Omega(\bold
x)=(n\phi_n(\frac1nx_n))_{n=1}^{\infty}$$
 (where we defined $\phi_n(x)=0$
if $\|x\|_{X_n}>1.$). We may argue that $\Omega$ is well-defined since
$$\frac1n\|x_n-u_n\|_{X_n}-\frac2{n^2}\le
\|\phi_n(\frac1nx_n)-\phi_n(\frac1nu_n)\|_{Y_n}\le
\frac1n\|x_n-u_n\|_{X_n}+\frac2{n^2}$$ if $\|x_n\|_{X_n},\|u_n\|_{X_n}\le
n.$  The same inequality implies that $\Omega$ is an isometry.

If $\bold y=(y_n)\in\ell_{\infty}(Y_n)$ we can define $x_n=n\psi_n(y_n/n)$
for all but finitely many $n.$  Then
$\|\phi_n(x_n/n)- y_n/n\|_{Y_n}\le \frac2{n^2}.$  It follows that
$\Omega$ is onto.

We also note that since $\|\phi_n(0)\|_{Y_n}\le \frac2{n^2}$ we have
$\Omega(0)=0.$  It follows from the Mazur-Ulam theorem (\cite{21}, \cite
{29}) that
$\Omega$ is linear.

Now notice that if $\|x_n\|_{X_n}=1$ then $\Omega(\bold
x)=(\Phi_n(x_n))$ (as elements of the ultraproduct).
Since $\Phi_n$ is homogeneous this means that for every $\bold x\in
X_{\Cal U}$ we have $\Omega(\bold x)=(\Phi_n(x_n)).$

Since (1) or (2) fails for every $n\ge 20$ we can assume that for some
set $P
\in\Cal U$ we have (1) failing for every $n\in P$ or (2) failing for
every $n\in P.$  If (1) fails there exists $\bold y=(y_n)$ with
$y_n\in S_{Y_n}$ such that if
$n\in P$ then $\|y_n-\Phi_n(x_n)\|\ge \epsilon$ whenever
$\|x_n\|_{X_n}\le 1.$  This contradicts the fact that $\Omega$ is an onto
isometry.

Similarly if (2) fails then there exist $\bold x_k=(x_{kn})$ for
$k=1,2,3$ with $\sum_{k=1}^3\|x_{kn}\|_{X_n}=1$, $\sum_{k=1}^3x_{kn}=0$
and $\|\sum_{k=1}^3\Phi_n(x_{kn})\|_{Y_n}\ge \epsilon$ for $n\in P.$
This contradicts the linearity of $\Omega.$\qed\enddemo

It is perhaps worth recording a result implicit in this argument:

\proclaim{Proposition 3.5}Suppose $(X_n)$ and $(Y_n)$ are two sequences
of real Banach spaces for which $\lim_{n\to\infty}d_{GH}(X_n,Y_n)=0.$
Then for any non-principal ultrafilter $\Cal U$ on $\Bbb N$ the
ultraproducts $\prod_{\Cal U}(X_n)$ and $\prod_{\Cal U}(Y_n)$ are
isometric.\endproclaim

Now we are able to establish that for $0<r<1$ we have equivalence for
the metrics $d_{GH}$ and $d_r$.

\proclaim{Theorem 3.6}Suppose $0<r<1$.  Then there
is a function $f=f_r:(0,1]\to (0,1]$ with $\lim_{\epsilon\to
0}f(\epsilon)=0$ such that for every pair of real Banach spaces $X$ and
$Y$ we have $$ 2^{-2/r+1}d_{GH}(X,Y) \le d_r(X,Y)\le f(d_{GH}(X,Y)).$$
\endproclaim

\demo{Remark}It would be interesting to find an explicit function $f$
satisfying the conditions of the theorem.\enddemo

\demo{Proof}In fact the upper estimate follows immediately from
Proposition 3.3 and Lemma 3.4.

For the lower estimate we begin by noticing that if $(\Phi,\Psi)\in\Cal
F_h(X,Y)$ and $\Delta_r(X,Y)=\sigma$ then for $x_1,\ldots,x_m\in X$ and
$y_1,\ldots,y_n\in Y$ we have, by putting
$u=\sum_{i=1}^mx_i-\sum_{j=1}^n\Psi(y_j),$
$$ \|\sum_{i=1}^m\Phi(x_i)-\sum_{j=1}^ny_j-\Phi(u)\|_Y \le 2^{1/r}\sigma
(\sum_{i=1}^m\|x_i\|_X^r+\sum_{j=1}^n\|y_j\|_Y^r)^{1/r}.$$

Since $\|\Phi(u)\|_Y \le \|u\|_X$ this implies that
$$ \|\sum_{i=1}^m\Phi(x_i)-\sum_{j=1}^ny_j\|_Y \le \|u\|_X+2^{1/r}\sigma
(\sum_{i=1}^m\|x_i\|_X^r+\sum_{j=1}^n\|y_j\|_Y^r)^{1/r}.$$

By using this and the symmetrical calculation with $X,Y$ interchanged we
obtain
$$\left
|\|\sum_{i=1}^mx_i-\sum_{j=1}^n\Psi(y_j)\|_X-\|\sum_{i=1}^m\Phi(x_i)
-\sum_{j=1}^ny_j\|_Y \right|\le
2^{1/r} \sigma
\left(\sum_{i=1}^m\|x_i\|^r_X+\sum_{j=1}^n\|y_j\|^r_Y\right)^{1/r}
.$$
If we let $\phi$ and $\psi$ be the restrictions of $\Phi$ and $\Psi$ to
the respective unit balls and apply the above inequality for $m+n\le 2$
we obtain that $D(\phi,\psi)\le 2^{2/r-1}\sigma.$
\qed\enddemo

We now recall the definition of $\Cal K-$space from \cite{16} or
\cite{18}.
 We say that a Banach space $X$ is a $\Cal K-${\it space} if there
is a constant $\kappa$ (we denote the best such constant by $\kappa
(X)$) such that whenever
$f:X\to\Bbb R$ is a homogeneous function satisfying
$$ |f(x_1+x_2)-f(x_1)-f(x_2)|\le \|x_1\|_X+\|x_2\|_X$$ there is a linear
functional $g:X\to\Bbb R$ with $|f(x)-g(x)|\le \kappa \|x\|_X$ for all
$x\in X.$
It is also natural to consider the notion of a $\Cal K_0-$space.
 We say that a Banach space $X$ is a $\Cal K_0-${\it space} if there
is a constant $\kappa_0$ (we denote the best such constant by $\kappa_0
(X)$) such that whenever
$f:X\to\Bbb R$ is a homogeneous function, which is bounded on $B_X$ and
satisfies
$$ |f(x_1+x_2)-f(x_1)-f(x_2)|\le \|x_1\|_X+\|x_2\|_X$$ there is a linear
functional $x^*\in X^*$ with $|f(x)-x^*(x)|\le \kappa_0 \|x\|_X$
for all
$x\in X.$  Clearly $\kappa_0(X)\le \kappa(X).$
If $X$ has the Bounded Approximation Property it may be shown that
$X$ is a $\Cal K_0$-space if and only if $X$ is a $\Cal K$-space. In
general however this equivalence is not known.

It is known that a Banach space $X$ is a $\Cal K-$space if $X$ has
nontrivial type (\cite{16}) or if $X^*$ is isomorphic to a subspace
of $L_1$ (see \cite{19}).

\proclaim{Theorem 3.7}Suppose $0<r<1$.  Then there is a constant $C=C(r)$
so that if $X$ is a $\Cal K_0-$space then for any Banach
space $Y$, $d_K(X,Y)\le C\kappa_0(X)d_r(X,Y).$

In particular if $X$ is a real $\Cal K_0-$space then for any sequence of
Banach spaces $(X_n)$ we have
$\lim_{n\to\infty}d_{GH}(X_n,X)=0$ implies
$\lim_{n\to\infty}d_K(X_n,X)=0.$\endproclaim

\demo{Remark}Notice that this theorem applies if $X$ is super-reflexive
or if $X$ is isomorphic to $c_0$ or $C(K)$ for some compact Hausdorff
space. As we note later in the case when $X$ is isomorphic to
$c_0$ or $\ell_{\infty}$ one can show that
$\lim_{n\to\infty}d_{BM}(X_n,X)=0$ whenever
$\lim_{n\to\infty}d_{GH}(X_n,X)=0$ (again for real spaces only).\enddemo

\demo{Proof}Suppose
$(\Phi,\Psi)\in \Cal F_h(X,Y),$ and $\sigma=\Delta_r(\Phi,\Psi).$ Then
for
$x_1,x_2\in X$ we have
$$ \|\Phi(x_1+x_2)-\Phi(x_1)-\Phi(x_2)\|_Y \le
2^{2/r-1}\sigma(\|x_1\|_X+\|x_2\|_X).$$
If $\|y^*\|_{Y^*}=1$ then there is a linear map $g:X\to\Bbb R$ such
that
$$ |y^*(\Phi(x))-g(x)|\le 2^{2/r-1}\sigma\kappa_0(X)\|x\|_X.$$  Now if
$\sum_{i=1}^nx_i=0$ we obtain
$$ |\sum_{i=1}^ny^*(\Phi(x_i))|\le
2^{2/r-1}\sigma\kappa_0(X)\sum_{i=1}^n\|x_i\|_X.$$  This in turn
implies that
$$ \|\sum_{i=1}^n\Phi(x_i)\|\le
2^{2/r-1}\sigma\kappa_0(X)\sum_{i=1}^n\|x_i\|_X.$$
 The result follows on appealing to Proposition 3.3. \qed\enddemo

\demo{Example}We now show by example that the Kadets and
Gromov-Hausdorff distances are not equivalent for general Banach spaces.
To see this we show that $\lim_{p\to 1} d_{GH}(\ell_p,\ell_1)=0.$
It is known (cf\cite{22} or \cite{24}) that $d_K(\ell_1,\ell_p)=1$ for
all $p>1.$

We consider the Mazur map $\phi:B_{\ell_p}\to B_{\ell_1}$ which is
defined by
$$ \phi(\xi)=(\sgn \xi_n |\xi_n|^p)_{n=1}^{\infty}$$
for $\xi=(\xi_n).$  Let $\psi=\phi^{-1}$ and then we compute
$D(\phi,\psi)$.  We need only estimate
$$ \left|\|\xi-\eta\|_{\ell_p}-\|\phi(\xi)-\phi(\eta)\|_{\ell_1}\right|$$
over all $\xi,\eta\in B_{\ell_p}.$
Clearly
$$ \left|\|\xi-\eta\|_{\ell_p}-\|\xi-\eta\|_{\ell_p}^p\right|\le 2^p-2.$$

 For any $a,b\in \Bbb R$ we have the estimate
$$ \left| |a-b|^p-|\sgn a|a|^p-\sgn b|b|^p|\right|\le (2^{p-1}-1)(|a|^p +
|b|^p).$$
To see this note that if $s,t\ge 0$ we have
$$ s^p+t^p\le (s+t)^p\le 2^{p-1}(s^p+t^p).$$
This implies that if $0\le s\le t$ then
$$ (t-s)^p \le t^p-s^p \le (t-s)^p+(2^{p-1}-1)(s^p+(t-s)^p).$$
The required inequality now follows by considering cases.

Now by summing we obtain
$$ \left|\|\xi-\eta\|_{\ell_p}^p-\|\phi(\xi)-\phi(\eta)\|_{\ell_1}\right|
\le 2^p-2.$$

This implies that $D(\phi,\psi)\le 2(2^p-2).$  Hence
$d_{GH}(\ell_1,\ell_p)\le 2(2^p-2)\to 0$ as $p\to 1.$\qed\enddemo

\heading{4. Continuity of certain maps for the Kadets
 metric} \endheading

In this section, we will establish a number of continuity-type results
for the Kadets distance.  Thus for example Theorem 4.3 below can be
interpreted as saying that the map $X\to X^*$ is continuous for the
Kadets distance.

\proclaim{Theorem 4.1} Let $Z$ be a Banach space and let $E,F$ be closed
subspaces of $Z$.  Then
$$ d_K(Z/E,Z/F) \le 2\Lambda(E,F).$$
\endproclaim

\demo{Remark}A somewhat different version of this result can be found in
the unpublished manuscript \cite{14}.\enddemo

\demo{Proof}Let $X=Z/E$ and $Y=Z/F.$ Suppose $\sigma>\Lambda(E,F).$ Let
$q_X$ and
$q_Y$ be the
respective quotient maps.  If $\theta>1$ we can define a homogeneous map
$f_X:X\to Z$ such that $q_Xf_X=I_X$ and $\|f_Xx\|_Z \le \theta\|x\|_X.$
We make a similar definition of $f_Y.$  Finally define $\Phi:X\to Y$ by
$\Phi(x)= \theta^{-1}q_Yf_X(x)$ and $\Psi:Y\to X$ by
$\Psi(y)=\theta^{-1}q_Xf_Y(y).$

Suppose $x_1,\ldots,x_m\in X$ and $y_1,\ldots,y_n\in Y.$  Let
$u=\sum_{i=1}^mx_i +\sum_{j=1}^n\Psi(y_j).$  Then
$f_X(u)-\sum_{i=1}^mf_X(x_i)-\theta^{-1}\sum_{j=1}^nf_Y(y_j)\in E.$  It
follows that there exists $z\in F$ with
$$ \|f_X(u)-\sum_{i=1}^mf_X(x_i)-\theta^{-1}\sum_{j=1}^nf_Y(y_j)-z\|_Z
\le \sigma(\theta(\|u\|_X+\sum_{i=1}^m\|x_i\|_X) + \sum_{j=1}^n\|y_j\|_Y.
)$$
Since $\|u\|_X \le \sum_{i=1}^m\|x_i\|_X +\sum_{j=1}^n\|y_j\|_Y$ this
implies that
$$ \|\theta\Phi(u)- \theta \sum_{i=1}^m\Phi(x_i) -\theta^{-1}\sum_{j=1}^n
y_j\|_Y \le 2\sigma \theta
(\sum_{i=1}^m\|x_i\|_X+\sum_{j=1}^n\|y_j\|_Y).$$

Note $\|\Phi(u)\|_Y\le \|u\|_X$ so that if
$v=\sum_{i=1}^m\Phi(x_i)+\sum_{j=1}^ny_j$ we obtain that
$$ \|v\|_Y \le \|u\|_X +  (2\sigma +
(1-\theta^{-2}))(\sum_{i=1}^m\|x_i\|_X+\sum_{j=1}^n\|y_j\|_Y).$$
Combined with the corresponding converse inequality this implies that
$$ \Delta(\Phi,\Psi)\le 2\sigma + (1-\theta^{-2}).$$
The result now follows.\qed\enddemo

It is interesting to note that there is a converse to this result.

\proclaim{Theorem 4.2}Suppose $X$ and $Y$ are Banach spaces, and suppose
$\sigma>d_K(X,Y).$
 Then
there is a Banach space $Z$ with closed subspaces $E,F$ such that
$Z/E$ is isometric to $X$, $Z/F$ is isometric to $Y$ and $\Lambda(E,F)\le
\sigma.$ Furthermore $E$ is isometric to $Y$ and $F$ is isometric to $X$.
\endproclaim

\demo{Proof}Suppose that $(\Phi,\Psi)\in \Cal
F_h(X,Y)$ with $\Delta(\Phi,\Psi)\le\sigma.$  We apply the construction
of Theorem 2.3.
As there, we
define $Z$ to be algebraically $X\oplus Y$ and then define the norm by
$$ \|(u,v)\|_Z=\inf\left(\|x_0\|_X+\|y_0\|_Y +\sigma
\sum_{i=1}^m\|x_i\|_X+\sigma\sum_{j=1}^n\|y_j\|_Y\right)$$ where the
infimum
is taken over all $x_0,\ldots,x_m\in X$ and $y_0,\ldots,y_n\in Y$ such
that $u=x_0+\sum_{j=1}^mx_j+\sum_{j=1}^n\Psi(y_j)$ and
$v=y_0+\sum_{j=1}^m\Phi(x_j)+\sum_{j=1}^ny_j.$

Now we have $\Lambda(X,Y)\le \sigma, $ where $X$ is
identified with $X\oplus\{0\}$ and
$Y$ with $\{0\}\oplus Y.$ Let us compute
$Z/Y:$ this is easily seen to be isometric to $X$ since
$$ \inf_{y\in Y}\|(x,y)\|_Z =\|(x,\Phi(x))\|_Z =\sigma\|x\|_X.$$
Similarly $Z/X$ is isometric to $Y$.
\qed\enddemo

\demo{Remark}If $X$ and $Y$ are separable then so is $Z$ and there is a
quotient map of $\ell_1$ onto $Z$.  It may then be verified that if
$\sigma>d_K(X,Y),$
$X$ and
$Y$ can be represented as $\ell_1/E$ and $\ell_1/F$  where
$\Lambda(E,F)\le 2\sigma.$   \enddemo

\proclaim{Theorem 4.3}Suppose $X$ and $Y$ are Banach spaces.  Then
$d_K(X^*,Y^*)\le 2d_K(X,Y).$\endproclaim

\demo{Remark} The reverse inequality (i.e.
$d_K(X,Y)\le cd_K(X^*,Y^*))$ is not valid.
In order to see this consider the case $X^*=Y^*=l_1,\
X=c_0,\ Y=$predual of $l_1$, non-isomorphic to $c_0$.
Proposition 5.3 of \cite{24} imply that $d_K(X,Y)\ge 1/2$.
On the other hand $d_K(X^*,Y^*)=0.$
\enddemo

\demo{Proof}Suppose $\sigma>d_K(X,Y).$  Then as in  Theorem 4.2, there is
a Banach space $Z$ and closed subspaces $E,F$ with $\Lambda(E,F)\le
\sigma$ and $Z/E=X,\ Z/F=Y.$  We then have $d_K(X^*,Y^*)\le
\Lambda(E^{\perp},F^{\perp}).$  We now use the fact that
$\Lambda(E^{\perp},F^{\perp})\le 2\Lambda(E,F)$ which is well-known (cf.
Theorem 3.4(d) and 3.13 of \cite{24}); we will provide a direct
proof. If
$z^*\in E^{\perp}$ and
$\|z^*\|_{Z^*}=1$ then for $f\in F$ we have $|z^*(f)|\le \sigma\|f\|_Z.$
By the Hahn-Banach theorem there exists $u^*\in Z^*$ with
$\|u^*\|_{Z^*}\le \sigma$ and $z^*-u^*\in F^{\perp}.$  Then
$\|z^*-u^*\|_{Z^*}\le 1+\sigma$ and so $d(z^*,B_{F^{\perp}})\le 2\sigma.$
It follows by symmetry that $\Lambda(E^{\perp},F^{\perp})\le
2\sigma.$\qed\enddemo

\proclaim{Theorem 4.4} If $X,Y$ are Banach spaces then we have \newline
(1) $d_K(X^{**},Y^{**}) \le d_K(X,Y)$\newline (2) For any fixed
ultrafilter $\Cal U$ on $\Bbb N$ $d_K(X_{\Cal U},Y_{\Cal U})\le d_K(X,Y)$
\newline (3)
$d_K(X^{**}/X,Y^{**}/Y)\le 8d_K(X,Y).$\endproclaim

\demo{Proof}For (1) observe simply that $X$ and $Y$ are embedded into a
common Banach space $Z$ then $X^{**},Y^{**}$ can be identified with
$X^{\perp\perp}, Y^{\perp\perp}$ in $Z^{**}$ and
$\Lambda(X^{\perp\perp},Y^{\perp\perp})\le \Lambda(X,Y).$  For (2) we may
use a similar argument with $Z_{\Cal U}.$
 For (3) we use
4.2. If
$\sigma>d_K(X,Y)$ there is a Banach space
$Z$
with closed subspaces $E,F$ with $\Lambda(E,F)\le \sigma$ and such that
$Z/E$ is isometric to $X$ and $Z/F$ is isometric to $Y.$ Furthermore from
the construction we have $Z=E+F.$ Consider
$Z^{**}$; then it is easy to verify that
$\Lambda(E^{\perp\perp},F^{\perp\perp})\le \sigma$.  Let $Q_X:Z\to
X$ and $Q_Y:Z\to Y$  be the quotient maps. Now suppose
$z^{**}\in Z+F^{\perp\perp}$ with $\|z^{**}\|_{Z^{**}}=1$ and that
$\epsilon>0.$ Then
$Q_Y^{**}z^{**}\in X\subset X^{**}$
so that there exists $z\in Z$ with $\|z\|_Z\le (1+\epsilon)$ and
$Q_Yz=Q_Y^{**}z^{**}$.  Thus $z^{**}-z\in F^{\perp\perp}.$  Now pick
$e^{**}\in E^{\perp\perp}$ with $\|z^{**}-z-e^{**}\|_{Z^{**}}\le
(2+\epsilon)\sigma.$  Thus $d(z^{**}, E^{\perp\perp}+Z)\le 2\sigma.$
This and the similar inequality with $E,F$ reversed leads to the estimate
$\Lambda(Z+F^{\perp\perp},E^{\perp\perp}+Z)\le 4\sigma$.  Now note that
$Z^{**}/(Z+F^{\perp\perp})$ is isometric to $Y^{**}/Y$ and
$Z^{**}/(Z+E^{\perp\perp})$ is isometric to $X^{**}/X.$ Now
Theorem
4.1 gives the result.
\qed\enddemo

We now turn to complex interpolation of Banach spaces. We
present an approach which encompasses several such situations.

Let $E$ be a complex Banach space and $\Cal U$ is an open subset
of the
complex plane $\Bbb C$ which is conformally equivalent to the unit disk
$\Cal D.$  We define an {\it interpolation field} $\Cal X$ to be a vector
space
of $E$-valued analytic functions $f:\Cal U\to E$ equipped with a norm
$\|\cdot\|_{\Cal X}$ such that $\Cal X$ is a Banach space and such that
the following two conditions hold:
\item{(1)} If $\varphi:\Cal U\to \Cal D$ is a conformal equivalence then we
have $f\in \Cal X$ if and only if $\varphi f\in \Cal X$ and $\|f\|_{\Cal
X}=\|\varphi f\|_{\Cal X}.$
\item{(2)} Each evaluation $f\to f(w)$ for $w\in\Cal U$ is bounded from
$\Cal X$ to $E.$

If $\Cal X$ is an interpolation field we define for each $w\in\Cal U$
the Banach space $X_w$ to be the set of $x\in E$ such that there exists
$f\in\Cal X$ with $f(w)=x$ equipped with the quotient norm
$\|x\|_w=\inf\{\|f\|_{\Cal X}:f(w)=x\}.$

This definition is easily seen to encompass the standard definitions of
complex interpolation spaces in the literature by taking $\Cal U$ to be
the open strip $\Cal S=\{w:0<\Re w<1\}$.  It also covers the
interpolation
method introduced by Coifman, Cwikel, Rochberg, Sagher and Weiss
in \cite{7}.

Our main result is that the map $w\to X_w$ is continuous for the Kadets
metric.   To state our result we introduce the pseudo-hyperbolic metric
on $\Cal U.$
 If $\xi,\eta\in \Cal U$ let $h(\xi,\eta)$
be the pseudo-hyperbolic distance on $\Cal U$ defined by
$h(\xi,\eta)=|\varphi(
\eta)|(=h(\eta,\xi))$ where $\varphi:\Cal
U\to\Cal D$ is a conformal equivalence with $\varphi(\xi)=0.$

\proclaim{Theorem 4.5}Let $\Cal U$ be an open subset of the complex
plane which is conformally equivalent to the unit disk, and let $\Cal F$
be an interpolation field on $\Cal U.$
  Then if $\xi,\eta\in\Cal U,$
$$ d_{K}(X_{\xi},X_{\eta})\le 2h(\xi,\eta).$$ \endproclaim

\demo{Proof}The argument is very similar to that of Theorem 4.1.
Suppose $\theta>1$.
 We define a homogenous map $\phi_{\xi}:X_{\xi}\to \Cal X$
with the property that $\phi_{\xi}(x;\xi)=x$ and $\|\phi_{\xi}(x)\|_{\Cal
X}\le
\theta\|x\|_{\xi}.$  Define $\Phi:X_{\xi}\to X_{\eta}$ by
$\Phi(x)=\theta^{-1}\phi_{\xi}(x;\eta)$
and similarly let $\Psi:X_{\eta}\to X_{\xi}$ be defined by
$\Psi(y)=\theta^{-1}\phi_{\eta}(y;\xi).$

Now suppose $x_1,\ldots,x_m\in X_{\xi}$ and $y_1,\ldots, y_n\in
X_{\eta}.$  Let $u=\sum_{i=1}^mx_i +\sum_{j=1}^n\Psi(y_j)$ and
$v=\sum_{i=1}^m\Phi(x_i)+\sum_{j=1}^ny_j.$   We note that
$$
\phi_{\xi}(u;\xi)-\sum_{i=1}^m\phi_{\xi}(x_i;\xi)-\theta^{-1}
\sum_{j=1}^n\phi_{\eta}(y_j;\xi)=0.$$

Let $\varphi$ be a conformal map of $\Cal U$ onto $\Cal D$ with
$\varphi(\xi)=0.$  Then there exists $f\in\Cal X$ with
$$\phi_{\xi}(u)-\sum_{i=1}^m\phi_{\xi}(x_i)-\theta^{-1}\sum_{j=1}^n
\phi_{\eta}(y_j)=\varphi f.$$
Using the  fact that $\|f\|_{\Cal X}=\|\varphi f\|_{\Cal X},$ we obtain
$$ \|\theta \Phi(u) -\theta\sum_{i=1}^m\Phi(x_i)-\theta^{-1}\sum_{j=1}^n
y_j\|_{\eta} \le 2\theta|\varphi(\eta)|
(\sum_{i=1}^m\|x_i\|_{\xi}+\sum_{j=1}^n\|y_j\|_{\eta}).$$
It follows that
$$ \|v\|_{\eta}\le \|u\|_{\xi} + (2h(\xi,\eta)+(1-\theta^{-2}))(
\sum_{i=1}^m\|x_i\|_{\xi}+\sum_{j=1}^n\|y_j\|_{\eta}).$$
With the symmetrical inequality this leads to
$$ \Delta(\Phi,\Psi)\le
2 h(\xi,\eta)+(1-\theta^{-2})$$
and the result follows.\qed\enddemo

We can apply these results to standard interpolation couples.  Let $X_0
,X_1$ be a Banach couple and let $X_{\theta}=[X_0,X_1]_{\theta}$ be the
standard complex interpolation space obtained by the Calderon method.

\proclaim{Corollary 4.6}If $0<\theta<\phi<1$ then
$$d_{K}(X_{\theta},X_{\phi})\le 2\frac{\sin(\pi(\phi-\theta)/2)}
{\sin(\pi(\phi+\theta)/2)}.$$\endproclaim

\demo{Proof}Define
$$\varphi(z)=\frac{\sin(\pi(z-\theta)/2)}{\sin(\pi(z+\theta)/2)}.$$  Then
$\varphi$ is a conformal mapping of the strip $\{z:0<\Re z<1\}$ onto the
open unit disk with $\varphi(\theta)=0.$\qed\enddemo

\proclaim{Corollary 4.7} If $1<p<q<\infty$ then for the (complex) spaces
$\ell_p,\ell_q$ we have:
$$ d_{K}(\ell_p,\ell_q)\le 2\frac{\sin(\pi(1/p-1/q)/2)}{\sin
(\pi(1/p +1/q)/2)}.$$\endproclaim

\demo{Remark}
In fact Corollary 4.7 holds for the corresponding real spaces.
One way to see this is to note that the pair $(\Phi,\Psi)$ of Theorem 4.4 can be
can be constructed even in the case $\theta=1$ and then map real
sequences to real sequences.\qed\enddemo

\demo{Remark}The estimate of Corollary 4.6 improves on previous
estimates
(\cite {15}, \cite {22}).  The best known lower estimate (\cite {23})
is $$d_K(\ell_p,\ell_q)\ge 2^{1/p-1}-2^{1/q-1}.$$\enddemo

\demo{Remark}We can also consider complex interpolation of quasi-Banach
spaces.  If we fix $0<r<1$ and allow $E$ in our definition to be an
$r$-Banach space it is easy to show that the map $\xi\to X_{\xi}$ is
continuous for the pseudo-metric $d_r.$  By interpolating between
$\ell_r$ when $r<1$ and $\ell_2$ one can then see that $\lim_{p\to
1}d_{GH}(\ell_p,\ell_1)=0$ as verified directly in the previous
section.\enddemo

\heading{The Kadets and Gromov-Hausdorff topologies}\endheading

Let $\aleph$ be any arbitrary but fixed cardinal.  Then we may consider
the pseudo-metric space $\Cal B_{\aleph}$ of all Banach spaces with
density character at most $\aleph$ with the Kadets or Gromov-Hausdorff
pseudo-metrics.  We note that there are examples of non-isomorphic Banach
spaces for which $d_K(X,Y)=0$.  An easy way to construct examples is to
fix $1<s<\infty$ and take a sequence $1<p_n<\infty$ with
$\lim_{n\to\infty}p_n=s,$ but $p_n\neq s$ for all $n.$  Then
$X=\ell_s\oplus_2\ell_2(\ell_{p_n})$ and
$Y=\ell_2(\ell_{p_n})$ satisfy $d_K(X,Y)=0.$  This follows easily from
the estimate from Corollary 4.6,
$$ d_K(\ell_2(\ell_{p_n}),\ell_2(\ell_{r_n}))\le 2\sup_{n}
\frac{|\sin(\pi(p_n^{-1}-r_n^{-1})/2)|}{\sin(\pi
(p_n^{-1}+r_n^{-1})/2)}.$$  Clearly these spaces are non-isomorphic.  If
one takes $s=1$ then one gets an example where $d_{GH}(X,Y)=0$ but
$d_K(X,Y)=1$ since $\ell_1$ embeds in $X$ but not $Y.$

We will be interested in this section in the topology of the
(pseudo-)metric spaces $\Cal B_{\aleph}$ with the Kadets or
Gromov-Hausdorff distances (we will sometimes use the
term metric
with the understanding that the spaces actually considered are the
Hausdorff quotients).  We consider a fixed cardinal $\aleph$ to avoid
certain set-theoretic problems; the collection of all Banach spaces
fails to be a set.  The most interesting choice of $\aleph$ is of course
$\aleph_0$ the set of separable Banach spaces; however in computing
duals, biduals etc. it is necessary to consider larger cardinals. We
observe that each set
$\Cal
B_{\aleph}$ is clopen (closed and open) in any larger $\Cal B_{\aleph'}
.$

Let $\Cal P$ be a property of Banach spaces; then for each cardinal
$\aleph$ we may consider the set $\Cal P=\Cal P_{\aleph}$ of all $X\in
\Cal B_{\aleph}$ with property $\Cal P$, so that we can think of $\Cal
P$ as a set.  We will say
(following
\cite{24}) that $\Cal P$ is {\it stable} if there exists a fixed
$\alpha>0$ so that $X \in\Cal P$ and $d_K(X,Y)<\alpha$ imply $Y\in\Cal
P.$   We will say that $\Cal P$ is respectively {\it clopen, open,
closed} if
(for each
$\aleph$) the set
$\Cal P$ is respectively clopen, open or closed for the Kadets
pseudo-metric. Obviously a stable
property is clopen; also the negation of a stable or clopen property is
also stable or clopen.

There are many known examples of stable properties.  Let us list
some:\newline
(1) $X$ is separable.\newline
(2) $X$ does not contain $\ell_1$.\newline
(3) $X$ is reflexive.\newline
(4) $X$ is super-reflexive.\newline
(5) $X$ has nontrivial type.\newline

See \cite {22}, \cite {23} and \cite {24}, where other stable
properties
are also discussed.
.
We do not know of any examples of clopen properties which are not stable.

The following Proposition is trivial from Theorem 4.3 and Theorem 4.4.
Notice that
this proposition was known for stable properties (cf. \cite {24},
\cite
{2}) and for open properties only under some restrictions (\cite
{24}, \cite{2}).

\proclaim{Proposition 5.1}Suppose $\Cal P$ is a stable (respectively
clopen,
open, closed) property; then the properties $\Cal P^*= \{X:X^*\in\Cal
P\}$, $\Cal P^{co}=\{X:X^{**}/X\in\Cal P\}$ and $\Cal
P^{\Cal U}=\{X:X_{\Cal U}\in\Cal P\}$ for some fixed ultrafilter $\Cal
U$ on
$\Bbb N$ are also stable
(respectively, clopen, open, closed)
.\endproclaim

Note that for (1) and (2) above this leads to new stable properties:

$\{X:X^* \text{ is separable}\}$ and $\{X:X^* \text{ does not contain }
\ell_1\}.$

Let us also mention some examples of open properties:\newline
(6) $X$ is isomorphic to $\ell_1.$ \newline
(7) $X$ is isomorphic to $\ell_{\infty}.$  \newline
(8) $X$ is injective.\newline
(9) $X$ is isomorphic to $c_0.$

We refer again to \cite{22}, \cite{23} and \cite {24}.  In fact, in
each
example it is easy to show additionally that the Kadets distance defines
a topology equivalent on the set to the Banach-Mazur distance, i.e.
$\lim_{n\to\infty}d_K(X_n,X)=0$ if and only if
$\lim_{n\to\infty}d_{BM}(X_n,X)=0.$
 Notice
that we also have that the following are open properties by Proposition
5.1:\newline
(10) $X^*$ is isomorphic to $\ell_1.$  \newline
(11) $X$ is a  $\Cal L_1$-space.

Let us add to the list the following simple further open properties:

\proclaim{Proposition 5.2}The properties ``$X$ is isomorphic to a
subspace of $c_0$'' and ``$X$ is isomorphic to a subspace of
$\ell_{\infty}$'' are open.\endproclaim

\demo{Proof}We prove only the former assertion.  Assume that $X$ is
a Banach space with Banach-Mazur distance less than $\lambda$ to  a
subspace of
$c_0.$ Suppose
$d_K(X,Y)<\frac{1}{2\lambda+1}.$  Then $Y$ is separable.  There is a
separable Banach space
$Z$ which contains $X$ and $Y$ isometrically so that
$\sigma=\Lambda(X,Y)<
\frac{1}{2\lambda+1}.$  Let $T:X\to c_0$ be a linear operator with
$\|x\|_X\le \|Tx\|_{c_0}\le \lambda \|x\|_X.$   Then $T$ may be extended
to an operator $T:Z\to c_0$ with $\|T\|\le 2\lambda$ (cf \cite{30}).
If $y\in Y$ then
there exists $x\in X$ with $\|x\|_X\le \|y\|_Y$ and $\|x-y\|_Z\le
\sigma \|y\|_Y.$  Then $$\|Ty\|_{c_0}\ge
\|Tx\|_{c_0}-2\lambda\sigma\|y\|_Y  \ge
(1-\sigma-2\lambda\sigma)\|y\|_Y.$$  Thus $Y$ is also isomorphic to a
subspace of $c_0.$\qed\enddemo

\proclaim{Theorem 5.3}The property ``$\kappa_0(X)<\infty$''
is open.
\endproclaim

\demo{Proof}Suppose $X$ is a Banach space with
$\kappa_0=\kappa_0(X)<\infty.$  Suppose $d_K(X,Y)<\sigma,$ where
$(14+56\kappa_0)\sigma<1.$ Then
there exists a bijective norm-preserving homogeneous map $\Omega:X\to Y$
such that for $x_1,\ldots,x_n\in X$ we have
$$ \left|\|\sum_{k=1}^n\Omega(x_k)\|_Y-\|\sum_{k=1}^nx_k\|_X\right|\le
14\sigma\sum_{k=1}^n\|x_k\|_X.$$

Now suppose $f:Y\to \Bbb R$ is a homogeneous map which is bounded on
$B_Y$
and such that $|f(y_1+y_2)-f(y_1)-f(y_2)|\le \|y_1\|_Y+\|y_2\|_Y.$  Since
$f$ is bounded a weak$^*$ compactness argument shows the existence of a
best approximation $y^*\in Y^*$ so that
$\sup_{y\in B_Y}|f(y)-y^*(y)|=M$ is minimized.

Let $\varphi=f-y^*.$  Suppose $x_1,x_2\in X$.  Then
$\|\Omega(x_1+x_2)-\Omega(x_1)-\Omega(x_2)\|_Y \le
28\sigma(\|x_1\|_X+\|x_2\|_X).$
Hence
$$|\varphi(\Omega(x_1+x_2))-\varphi(\Omega(x_1)+\Omega(x_2))|\le
(28\sigma M+2) (\|x_1\|_X+\|x_2\|_X).$$
Now
$$
|\varphi(\Omega(x_1)+\Omega(x_2))-\varphi(\Omega(x_1))-
\varphi(\Omega(x_2))| \le  \|x_1\|_X+\|x_2\|_X$$  so that we conclude
that
$$
|\varphi(\Omega(x_1+x_2))-\varphi(\Omega(x_1))-\varphi(\Omega(x_2))|\le
(28\sigma M+3)(\|x_1\|_X+\|x_2\|_X).$$
It follows that there exists $x^*\in X^*$ so that
$$ |\varphi(\Omega(x))-x^*(x)| \le \kappa_0(28\sigma M +3)\|x\|_X.$$

Now suppose $y_1,\ldots,y_n\in Y$ and $\sum_{k=1}^ny_k=0.$  Let
$x=\sum_{k=1}^n\Omega^{-1}(y_k).$  Then $\|x\|_X \le
14\sigma \sum_{k=1}^n\|y_k\|_Y.$

Notice that
$$
\align
 |\sum_{k=1}^n
\varphi(y_k)-\varphi(\Omega(x))|&=|\sum_{k=1}^n\varphi
(\Omega\Omega^{-1}y_k) -\varphi(\Omega(x))-
x^*(\sum_{k=1}^n\Omega^{-1}(y_k)-x)| \\
&\le \kappa_0(3+28\sigma M)(1+14\sigma)\sum_{k=1}^n\|y_k\|_Y.
\endalign
$$

Now $$|\varphi(\Omega(x))|\le M\|x\|_X \le 14\sigma
M\sum_{k=1}^n\|y_k\|_Y
$$
so that we finally obtain
$$ |\sum_{k=1}^n\varphi(y_k)| \le C
\sum_{k=1}^n\|y_k\|_Y$$
where $$C=
\left(14\sigma
M+\kappa_0(6+56\sigma M)\right).$$

Now define a sublinear functional on $Y$ by
$$p(y)=\inf\left(\sum_{k=1}^n\varphi(y_k)+C\sum_{k=1}^n\|y_k\|_Y:\
\sum_{k=1}^ny_k=y\right).$$
Notice that if $\sum_{k=1}^ny_k=y$ then
$$ \sum_{k=1}^n\varphi(y_k)+C\sum_{k=1}^n\|y_k\|_Y \ge
\varphi(y)-C\|y\|_Y
$$ so that $p$ is well-defined.  Let $h$ be any linear functional on $Y$
so that $h(y)\le p(y)$ for all $y.$   Then $h(y)\le \varphi(y)+C\|y\|_Y$
and by applying to $-y$ we have $h(y)\ge \varphi(y)-C\|y\|_Y.$  Thus
$|h(y)-\varphi(y)|\le C\|y\|_Y.$

Now considering $f-y^*-h$ we see that we must have $M\le C$ and so if
$(14+ 56\kappa_0)\sigma<1$ we have
$$ M \le
\frac{6\kappa_0}{1-14\sigma-56\sigma\kappa_0}.$$
and this gives an estimate for $\kappa_0(Y).$\qed\enddemo

The notions of stable and open properties are most naturally applied to
interpolation scales.

\proclaim{Proposition 5.4}Suppose $(X_0,X_1)$ is a Banach couple and
that
$X_{\theta}=[X_0,X_1]_{\theta}.$  If $\Cal P$ is a stable property and
there exists $0<\theta<1$ so that $X_{\theta}$ has property $\Cal P$ then
$X_{\phi}$ has $\Cal P$ for every $0<\phi<1.$\endproclaim

\demo{Proof}This is immediate from Corollary 4.6. \qed\enddemo

\demo{Remark} This can be applied to each of the stable properties listed
above.  This yields a number of results, many of which are certainly
known to specialists.  For the example the case of reflexivity can be
deduced from Caldero\'n's original paper \cite{6}.  However we feel this
general framework for such a result has some interest.\enddemo

\proclaim{Proposition 5.5}
Suppose $(X_0,X_1)$ is a Banach couple and that
$X_{\theta}=[X_0,X_1]_{\theta}.$  If $\Cal P$ is an open property
and
there exists $0<\theta<1$ so that $X_{\theta}$
 has property $\Cal P$ then there exists $\epsilon>0$ so that
$X_{\phi}$ has $\Cal P$
 for every $\phi$ with
$|\phi-\theta|<\epsilon.$\endproclaim

Here we single out three special cases which seem to be new and of some
interest:

\proclaim{Proposition 5.6}
Suppose $(X_0,X_1)$ is a Banach couple and that
$X_{\theta}=[X_0,X_1]_{\theta}.$  Suppose
there exists $0<\theta<1$ so that
$X_{\theta}$ is isomorphic to $c_0$
(respectively isomorphic to a subspace of $c_0$, resp. isomorphic to
$\ell_1$)
   then there exists $\epsilon>0$ so that
$X_{\phi}$ is isomorphic to $c_0$
(respectively isomorphic to a subspace of $c_0$, resp. isomorphic to
$\ell_1$)
 for every $\phi$ with
$|\phi-\theta|<\epsilon.$\endproclaim

\demo{Remarks}Of course, Propositions 5.4-5.6 apply to general
interpolation fields.\enddemo

Proposition 5.4 suggests it is of interest to make the following
definition.  Let $X$ be an arbitrary Banach space.  We will say that the
{\it (Kadets) component, $\Cal C_X$} of $X$ is the intersection of all
clopen
properties containing $X.$  Clearly if $(X_0,X_1)$ is a Banach couple
then for all $0<\theta,\phi<1$ we have that $X_{\phi}$ is in the
component of $X_{\theta}.$

We first state some elementary properties of components.

\proclaim{Proposition 5.7}Let $X$ be an arbitrary Banach space.  Then:
\newline
(1) $Y\in\Cal C_X$ if and only if $\Cal C_Y=\Cal C_X.$\newline
(2) For every $Y\in\Cal C_X$ the density character dens $X=$ dens
$Y$ and dens $X^*=$ dens $Y^*.$\newline
(3) If $Y\in\Cal C_X$ and $Y_1$ is isomorphic to $Y$ then $Y_1\in\Cal
C_X.$
\newline
(4) If $Y$ is an arbitrary Banach space then for any subspace $E$ we have
$Y/E\oplus E\in\Cal C_X$ if and only if $Y\in\Cal C_X.$\endproclaim

\demo{Remark}Note that (3) allows us not to specify any special norm on
the spaces in (4).\enddemo

\demo{Proof}(1) is elementary.  For (2) note that the sets
$\{Y:\text{dens }Y=\text{dens }X\}$ and $\{Y: \text{dens }Y^*=\text{dens
} X^*\}$ are both clopen.

For (3) observe as in Proposition 6.7 of \cite{24} that there is a
family of isomorphic copies of $Y$, $Y_t$ for $0\le t\le 1,$ say, so that
$t\to Y_t$ is continuous for the Kadets distance, with $Y_0=Y.$

For (4) we use Lemma 5.9 of \cite{24}.  Let $Z=Y\oplus_1 Y/E$ and
let
$Q:Y\to Y/E$ be the quotient map.  For any $t\in\Bbb R$  with $t\neq
0$ we let
$G_t=\{(ty,Qy):y\in Y\}.$  Then it is easy to show that $\lim_{s\to
t}\Lambda (G_s,G_t)=0$ for any $t\in \Bbb R.$ If $t=0$ we define
$G_0=E\oplus_1 Y/E.$   Then using Lemma 5.9 of \cite{24} we also have
$\lim_{s\to 0}\Lambda(G_s,G_0)=0.$  Thus the map
$t\to G_t$
is continuous for the Kadets distance.  However $G_t$ is isomorphic
to
$Y$
for all $t\neq 0$ while $G_0=E\oplus_1 Y/E.$\qed\enddemo

\proclaim{Proposition 5.8}
We have $X\in \Cal C_{\ell_1}$ if and only if $X$ is separable and
contains a copy of $\ell_1.$\endproclaim

\demo{Proof}Since the set of $X$ which is both separable and contains a
copy of $\ell_1$ is stable inclusion is immediate.  Conversely suppose
$X$ is any separable Banach space containing a subspace $E$ isomorphic to
$\ell_1.$  Then $X$ is in the same component as $\ell_1\oplus_1 X/E.$
Let $F$ be the kernel of  a quotient map from $\ell_1$ onto $X/E.$
Let $G=\ell_1(F)\oplus \ell_1(X/E)\in\Cal C_{\ell_1}.$  Since the map
$Y\to Y\oplus_1 X/E$ is trivially seen to be continuous for the Kadets
metric the set $\{Y\oplus_1 X/E; Y\in\Cal C_{\ell_1}\}$ is connected and
meets $\Cal C_{\ell_1}$ since $G\oplus_1 (X/E)$ is isomorphic to $G.$
Thus $\ell_1\oplus_1 X/E \in \Cal C_{\ell_1}$ and so $X\in\Cal
C_{\ell_1}.$\qed\enddemo

\demo{Problem 1}What is the component of $c_0?$  Similar techiques to the
above Proposition show that this contains all infinite-dimensional
subspaces of $c_0.$  Proposition 5.7 (4) does not help to
give any other examples, since being a subspace of $c_0$ is a three-space
property (\cite{1}).\enddemo

\demo{Problem 2}What is the component of $\ell_2?$  This is contained in
the stable set of all separable super-reflexive spaces.  We do not
know
if it coincides with this set.  This is related to Pisier's notion of
$\theta-$Hilbertian spaces (\cite{28}); any space which is
$\theta-$Hilbertian for $\theta>0$ belongs to $\Cal C_{\ell_2}$.\enddemo

One can obviously introduce the notion of a Gromov-Hausdorff component in
an analogous fashion.  Clearly the Gromov-Hausdorff component of $X$,
$\Cal G_X$ contains $\Cal C_X.$  Since the map $p\to\ell_p$ is continuous
for the Gromov-Hausdorff distance for $p\in [1,\infty)$ we have
that
$\Cal
G_{\ell_1}=\Cal G_{\ell_2}.$  It seems quite possible that this will
correspond with the collection of all separable Banach spaces so we ask:

\demo{Problem 3}Is the set of all separable Banach spaces connected for
the Gromov-Hausdorff distance?
We are unable to decide if $c_0$ is in $\Cal G_{\ell_2}$ so it would
be interesting to identify $\Cal G_{c_0}.$
\enddemo

Let us take this opportunity to make a few remarks about the
Gromov-Hausdorff distance.

\proclaim{Proposition 5.9}The following properties are open for the
Gromov-Hausdorff distance:\newline
(1) $X$ is isomorphic to
$c_0$\newline (2) $X$ is isomorphic to $\ell_{\infty}$,
\newline
(3) $\kappa_0(X)<\infty$\newline
(4) $X$ has nontrivial type.
\newline
(5) $X$ has nontrivial cotype.\newline
(6) $X^*$ has nontrivial type.\newline
(7) $X^*$ has nontrivial cotype.

Furthermore, on the sets defined by properties (1) and (2)
Gromov-Hausdorff distance defines a topology equivalent to the
Banach-Mazur distance.
\endproclaim

\demo{Proof}(1) and (2)  follow almost immediately from
Theorem
3.7, and the remarks above. In a similar way (3) follows from Theorem 3.7
and Theorem 5.3. Then (4) follows since if $X$ has nontrivial type
then $\kappa_0(X)<\infty$ again by using Theorem 3.7.

To establish (5),(6) and (7) we use Proposition 3.5. Indeed it is
immediate from 3.5 that the property ``$\ell_{\infty}$ is finitely
representable in $X$'' is closed for the Gromov-Hausdorff distance.  Let
us establish (6) and (7) by showing that for a fixed Banach space $E$ the
condition ``$E$ is finitely representable in $X^*$'' is closed for
the Gromov-Hausdorff distance.  Then we may take $E=\ell_1$ and
$E=\ell_{\infty}.$

To show this last statement it suffices to suppose $E$
finite-dimensional. Then if $X_n$ converges to $X$ in
Gromov-Hausdorff distance i.e. $d_{GH}(X_n,X)\to
0$ and
$E$ is finitely representable in each $X_n^*$ it is immediate that $E$ is
isometric to a subspace of $(\prod_{\Cal U}(X_n))^*.$  Hence by Theorem
3.5 $E$ embeds isometrically into the space $X_{\Cal U}^*.$  Now
$(X^*)_{\Cal U}$ embeds naturally into $X_{\Cal U}^*$ as a norming
subspace by the identification
$$ \bold x^*(\bold x)=\lim_{n\in\Cal U}x_n^*(x_n).$$
The space of operators $\Cal L(E,X_{\Cal U}^*)$ is naturally a dual
space of $E\otimes_{\pi}X_{\Cal U}$ and $\Cal L(E, (X^*)_{\Cal U})$ is
norming as a subspace.  Hence $J:E\to X_{\Cal U}^*$ is an isometry we can
find a net $J_{\alpha}:E\to (X^*)_{\Cal U}$ so that $J_{\alpha}\to J$
weak$^*$ and $\|J_{\alpha}\|\le 1.$  By the weak$^*$ lower-semicontinuity
of the norm in $X_{\Cal U}^*$ it follows that for any $\epsilon>0$ there
exists
$\alpha$ with $(1-\epsilon)\|e\|_E\le \|J_{\alpha}\|_{(X^*)_{\Cal U}}\le
\|e\|_E$ for $e\in E.$   Hence $E$ is finitely representable in
$(X^*)_{\Cal U}$ and thus also in $X^*.$
\qed\enddemo

\demo{Remark}It has been conjectured by the first author that
$\kappa_0(X)<\infty$ might be equivalent to the property that $X^*$ has
finite cotype.
\enddemo

We conclude the paper by showing that the collection of all separable
Banach spaces is not itself separable for either the Gromov-Hausdorff or
Kadets distances.  More precisely we show:

\proclaim{Theorem 5.10}Suppose $1<p<\infty$ and $p\neq 2.$  Then the set
of Banach spaces $X$ isomorphic to $\ell_p$ is not separable in either
Kadets or Gromov-Hausdorff distances.\endproclaim

\demo{Proof}First notice that since $\kappa_0(\ell_p)<\infty$ the two
pseudo-metrics $d_K$ and $d_{GH}$ define equivalent topologies on the
set
$\Cal I_p$ of all isomorphic copies of $\ell_p.$  Secondly since $X\to
X^*$ is a homeomorphism on the set of reflexive spaces for the Kadets
pseudo-metric it suffices to consider $1<p<2.$

We shall consider spaces $\ell_p(E_n)$ where each $E_n$ is a
finite-dimensional Hilbert space.  Each such space is isomorphic to
$\ell_p$ by an old result of Pe\l czy\'nski \cite{26}.  We will prove
the
following Lemma:\enddemo

\proclaim{Lemma 5.11}There exists $\sigma_0>0$ so that if $(E_n)$ and
$(F_n)$ are two sequences of finite-dimensional Hilbert spaces and
$d_K(\ell_p(E_n),\ell_p(F_n))<\sigma_0$ then there is a bijection
$\pi:\Bbb N\to\Bbb N$ so that $\frac14\dim E_n\le \dim F_{\pi(n)}\le 10
\dim
E_n.$\endproclaim

\demo{Proof of the lemma}We will write $X=\ell_p(E_n)$ and
$Y=\ell_p(F_n).$ A typical element of $X$ will be denoted $\bold
x=(x_n)_{n=1}^{\infty}$ with a similar notation for $Y.$ We also adopt
the convention that
$\eta$ will denote a
function of $\sigma$ satisfying $\lim_{\sigma\to 0}\eta(\sigma)=0$ which
may vary from occurence to occurence.

First observe that by Theorem 2.4
there exists a homogeneous, norm-preserving bijection
$\Omega:X\to Y$ with the property that if $\bold x_1,\ldots,\bold x_n\in
X$ then $$ \left|\|\sum_{k=1}^n\Omega \bold x_k\|_Y-\|\sum_{k=1}^n\bold
x_k\|_X\right|\le 14\sigma\sum_{k=1}^n\|\bold x_k\|_X.$$

Suppose $\bold u,\bold v \in S_X$ have disjoint supports, i.e.
$\|u_n\|_{E_n}\|v_n\|_{E_n}=0$ for each $n.$  Then $\|\bold u\pm \bold
v\|_X=2^{1/p}$ and hence $\|\Omega(\bold u)\pm \Omega(\bold v)\|_Y \ge
2^{1/p}-28\sigma.$  Let $\Omega(\bold u)=\bold y$ and
$\Omega(\bold v)=\bold z.$ Then we have
$$
\align
2^{1/p}-28\sigma &\le (\frac12\|\bold y+\bold z\|_Y^2+\frac12\|\bold
y-\bold z\|_Y^2)^{1/2} \\
&=\left( (\sum_{k=1}^{\infty}2^{-p/2}\|y_k+z_k\|_{F_k}^p)^{2/p}
+(\sum_{k=1}^{\infty}2^{-p/2}\|y_k-z_k\|_{F_k}^p)^{2/p}\right)^{1/2}\\
&\le \left( \sum_{k=1}^{\infty}(\frac12 \|y_k+z_k\|_{F_k}^2
+\frac12\|y_k-z_k\|_{F_k}^2)^{p/2}\right)^{1/p}\\
&=
\left( \sum_{k=1}^{\infty}( \|y_k\|_{F_k}^2
+\|z_k\|_{F_k}^2)^{p/2}\right)^{1/p}\\
&\le
\left(\sum_{k=1}^{\infty}\max(\|y_k\|_{F_k},\|z_k\|_{F_k})^p
\right)^{1/p-1/2}
\left(\sum_{k=1}^{\infty}(\|y_k\|_{F_k}^p+\|z_k\|_{F_k}^p)\right)^{1/2}.
\endalign
$$

This implies an estimate that
$$
\left(\sum_{k=1}^{\infty}\max(\|y_k\|_{F_k},\|z_k\|_{F_k})^p
\right)^{1/p} \ge 2^{1/p}-\eta(\sigma)$$
and hence an estimate
$$ \left(\sum_{k=1}^{\infty}\min
(\|y_k\|_{F_k},\|z_k\|_{F_k})^p\right)^{1/p}\le \eta(\sigma).$$

Now suppose $\bold y\in S_Y$ is supported on exactly one co-ordinate $m$
say. Let $\Omega^{-1}(\bold y)=\bold x.$  Then we can write $\bold
x=\bold u+\bold v$ disjointly where
$\|\bold u\|_X\ge \|\bold v\|_X,$ and $\|\bold u\|_X^p\le \|\bold
v\|_X^p+\max\|x_n\|_{E_n}^p.$
Let $\bold w=\Omega(\bold u)$ and $\bold z=\Omega(\bold v).$  Then since
$\|\bold u\|_X,\|\bold v\|_X\le 1$ we obviously have that
$\min (\|w_m\|_{F_m},\|z_m\|_{F_m})\le \eta.$
Now
$$
\|\bold y- \Omega(\bold u)-\Omega(\bold v)\|_Y \le 42\sigma$$
so that
$$ \|y_m-w_m-z_m\|_{F_m}\le 42\sigma$$
and hence
$$ \|w_m+z_m\|_{F_m}\ge 1-42\sigma.$$
From this we have an estimate $\max(\|w_m\|_{F_m},\|z_m\|_{F_m})\ge
1-\eta.$  This in turn means $\|\bold u\|_X \ge 1-\eta$ and hence an
estimate
$\max\|x_n\|_{F_n}\ge 1-\eta.$

It follows that if $\sigma$ is small enough there is a unique $n=n(\bold
y)
$ so that
$\|x_n\|_{F_n}\ge 1-\eta.$

Now suppose that $\bold y_1$ and $\bold y_2$ are both unit vectors with
the same singleton support $m.$  If $n(\bold y_1)\neq n(\bold y_2)$ we
will have $\|\Omega^{-1}(\bold y_1)\pm \Omega^{-1} (\bold y_2)\|_X\ge
2^{1/p}-\eta$
but $\min(\|\bold y_1\pm \bold y_2\|_Y)\le 2^{1/2}.$  Again for
$\sigma$
small enough this is a contradiction.  Hence we conclude that $n$ is a
function only of $m.$

Applying the same technique to $X$ and $Y$ interchanged we obtain two
maps $\pi:\Bbb N\to\Bbb N$ and $\rho:\Bbb N\to\Bbb N$ such that if
$\bold y=\Omega\bold x$ and both are unit vectors then:\newline
(a) if $\bold x$ is supported only at $m$ then
$\|y_{\pi(m)}\|_{F_{\pi(m)}}\ge
1-\eta$ and\newline
(b) if $\bold y$ is supported only at $m$ then
$\|x_{\rho(m)}\|_{E_{\rho(m)}}\ge 1-\eta.$

Now if $\bold x$ is supported only at $m$ then $\bold y=\bold z+\bold w$
where $\|\bold z\|_Y\ge 1-\eta$ and $\|\bold w\|_Y\le \eta$ are disjoint
and $\bold z$ is supported only at $\pi(m).$  Thus
$$\|\bold x-\Omega^{-1}(\bold z)-\Omega^{-1}(\bold w)\|_X \le 42\sigma$$
from which it follows that for $\sigma$ small enough we must have
$\rho(\pi(m))=m$ and similarly $\pi(\rho(m))=m.$  Thus $\pi$ is a
bijection.

Pick a maximal subset $(\xi_k)_{k\in J}$ of $B_{E_m}$ with
$\|\xi_k-\xi_l\|_{E_m}> \frac12.$  Then $|J|\le 5^{\dim E_m}.$
Let $\bold x_k$ be the element of $X$ with zeros everywhere except
$\xi_k$ in the
$m$th position.   Suppose $\zeta\in S_{F_{\pi(m)}}$; let $\bold z$ be
the similarly defined element of $Y$ with exactly one nonzero element
$\zeta.$  Let $\bold x=\Omega^{-1}\bold z.$  Then there exists $k\in J$
with $\|\bold x-\bold x_k\|_X \le \frac12 +\eta(\sigma).$   Then
$$ \|\bold z -\Omega(\bold x_k)-\Omega(\bold x-\bold x_k)\|_Y \le
42\sigma.$$
Hence
$$ \|\bold z-\Omega(\bold x_k)\|_Y \le \frac12+\eta(\sigma).$$

Now let $(\psi_k)$ be the $\pi(m)-$co-ordinate of $\Omega(\bold x_k).$
We have
$$ \|\zeta-\psi_k\|_{F_{\pi(m)}} \le \frac12+\eta(\sigma).$$
For $\sigma$ small enough this implies that $(\psi_k)_{k\in J}$ is a
3/4-net in $B_{F_{\pi(m)}}$ so that $$|J|\ge (4/3)^{\dim F_{\pi(m)}}.$$  We
conclude that $\dim F_{\pi(m)}\le C\dim E_m$ where $C=\log 5/\log
(4/3)\le 10.$  This proves the lemma.
\qed\enddemo

\demo{Proof of the theorem}For each infinite subset $M$ of $\Bbb N$ form
the space $X_M=\ell_p(\ell_2^{10^{2n}})_{n\in M}.$  The lemma shows these
are uniformly separated.\qed\enddemo

\enddocument

\end